\documentclass[10pt, a4paper]{amsart}
\usepackage{amsthm}
\theoremstyle{}
{\theoremstyle{definition}
\newtheorem{dfn}{Definition}[section]}
\newtheorem{prop}[dfn]{Proposition}
\newtheorem{nota}[dfn]{Notation}
\newtheorem{thm}[dfn]{Theorem}
{\theoremstyle{definition}
\newtheorem{rem}[dfn]{Remark}}
\newtheorem{lem}[dfn]{Lemma}
\newtheorem{cor}[dfn]{Corollary}

{\theoremstyle{definition}
\newtheorem{exa}[dfn]{Example}}

\usepackage[top=30truemm,bottom=30truemm,left=30truemm,right=30truemm]{geometry}
\usepackage{amsmath,amssymb,amscd,bm,graphicx,tikz-cd,dsfont}
\usepackage[all]{xy}
\usepackage{hyperref}

\newcommand{\MF}{\operatorname{MF}}

\newcommand{\DMF}{\operatorname{DMF}}
\newcommand{\KMF}{\operatorname{KMF}}
\newcommand{\AMF}{\operatorname{AMF}}

\newcommand{\Dsg}{\operatorname{D^{sg}}}
\newcommand{\Sing}{\operatorname{Sing}}
\newcommand{\Spec}{\operatorname{Spec}}
\newcommand{\Supp}{\operatorname{Supp}}
\newcommand{\Hom}{\operatorname{Hom}}
\newcommand{\Max}{\operatorname{Max}}
\newcommand{\Min}{\operatorname{Min}}
\newcommand{\Singloc}{\operatorname{Sing^{loc}}}

\newcommand{\coh}{\operatorname{coh}}
\newcommand{\Qcoh}{\operatorname{Qcoh}}
\newcommand{\Mod}{\operatorname{Mod}}

\newcommand{\Perf}{\operatorname{Perf}}

\newcommand{\CM}{\operatorname{CM}}
\newcommand{\codim}{\operatorname{codim}}
\newcommand{\embdim}{\operatorname{emb.dim}}
\newcommand{\normalbar}{\operatorname{-}}

\newcommand{\cA}{\mathcal{A}}

\newcommand{\cK}{\mathcal{K}}
\newcommand{\cL}{\mathcal{L}}

\newcommand{\cO}{\mathcal{O}}

\newcommand{\cT}{\mathcal{T}}

\newcommand{\bC}{\mathbb{C}}

\newcommand{\m}{\mathfrak{m}}
\newcommand{\n}{\mathfrak{n}}

\newcommand{\colim}{\varinjlim}
\renewcommand{\lim}{\varprojlim}

\begin{document}

\title[]{Relative singular locus  and   Balmer spectrum of matrix factorizations}

\author[Y.~Hirano]{Yuki Hirano}
\date{}


\begin{abstract} 

 For a separated Noetherian scheme $X$ with an ample family of line bundles and a non-zero-divisor  $W\in\Gamma(X,L)$ of a line bundle $L$ on $X$, we classify certain  thick subcategories of the derived matrix factorization category ${\rm DMF}(X,L,W)$ of the Landau-Ginzburg model $(X,L,W)$. Furthermore, by using the classification result and the theory of Balmer's tensor triangular geometry,
we show that the spectrum of the tensor triangulated category $({\rm DMF}(X,L,W), \otimes^{\frac{1}{2}})$ is  homeomorphic to the relative singular locus $\Sing (X_0/X)$,  introduced in this paper, of the zero scheme $X_0\subset X$ of $W$.

\end{abstract}

\maketitle{}

\section{Introduction}

\subsection{Background}
For a given category of algebraic objects associated to a scheme, it is expected that we can  extract geometric information of the scheme or the scheme itself from the category.  Gabriel  reconstructed a Noetherian scheme $X$ from the abelian category $\Qcoh X$ of quasi-coherent sheaves on $X$  \cite{gab}, and later Rosenberg generalized the reconstruction theorem for arbitrary schemes \cite{rose}.  
Although we can't reconstruct a smooth variety  from the derived category of coherent sheaves in general,
Balmer  reconstructed arbitrary Noetherian scheme $X$ from the tensor triangulated category $(\Perf X,\otimes)$ of perfect complexes on $X$ with  the natural tensor structure $\otimes$ \cite{balmer}. Balmer's idea is to associate to any tensor triangulated category $(\cT,\otimes)$ a ringed space $\Spec (\cT,\otimes)=({\rm Spc}(\cT,\otimes),\cO)$, and he proved an isomorphism $X\cong \Spec(\Perf X,\otimes)$ by using Thomason's result of classification of thick subcategories  of perfect complexes $\Perf X$ which are closed under $\otimes$-action of $\Perf X$.

In addition to the Thomason's result, classifications of thick subcategories of triangulated categories are studied in many articles. 
For example, Takahashi classified thick subcategories of the stable category $\underline{\CM}(R)$ of maximal Cohen-Macaulay modules over  an abstract hypersurface local  ring $R$ \cite{takahashi}. Stevenson proved a classification  of certain thick subcategories of the singularity category $\Dsg(X)$ of a hypersurface singularity $X$ \cite{stev}.



\subsection{Relative singular locus}
To state our main results, we introduce a new notion of relative singular locus. Let $i:T\hookrightarrow S$ be a closed immersion of Noetherian schemes. We define the {\it relative singular locus}, denoted by  $\Sing(T/S)$, of $i$ as the following subset of $T$;
\[
\Sing(T/S):=\{\,p\in T\,\mid\, \exists F\in \coh T \hspace{2mm}{\rm such\hspace{1mm} that}\hspace{2mm} F_p\notin\Perf \cO_{T,p} \hspace{2mm}{\rm and}\hspace{2mm} i_*(F)\in \Perf S\,\}.
\] 
We also consider the {\it locally relative singular locus} $\Singloc(T/S)$ of $i$ defined by 
\[
\Singloc(T/S):=\{\,p\in T\,\mid\, \exists M\in{\rm mod}\,{\cO_{T,p}} \hspace{2mm}{\rm such\hspace{1mm} that}\hspace{2mm} M\notin\Perf \cO_{T,p} \hspace{2mm}{\rm and}\hspace{2mm} i_{p*}(M)\in \Perf \cO_{S,p}\,\}.
\] 
 By definition, we have the inclusions 
 \[\Sing(T/S)\subseteq\Singloc(T/S)\subseteq\Sing(T),\]
where $\Sing(T)$ is the usual singular locus of $T$. These loci can be different to each other, but,  if $S$ is regular,  these loci are equal; $\Sing(T/S)=\Singloc(T/S)=\Sing(T)$.   
Roughly speaking, the relative singular locus $\Sing(T/S)$ is a set of points $p$ in $T$ such that the mildness of the singularity of $p$ in $T$ is worse than the mildness of the singularity of $p$ in $S$. 
In fact, for a quasi-projective variety $X$ over $\bC$ and a regular function $f\in\Gamma(X,\cO_X)$ which is non-zero-divisor, we have the following equality of subsets of the associated complex analytic space $(X^{\rm an},\cO_{X}^{\rm an})$;
\[
\Singloc(f^{-1}(0)/X)\cap X^{\rm an}={\rm Crit}(f^{\rm an})\cap{\rm Zero}(f^{\rm an}),
\]
 where ${\rm Crit}(f^{\rm an})$ denotes the  critical locus of the associated  function $f^{\rm an}\in\Gamma(X^{\rm an},\cO_{X}^{\rm an})$, which is  defined by 
 \[
 {\rm Crit}(f^{\rm an}):=\{\, p\in X^{\rm an}\,|\, ({f^{\rm an}})_p-f^{\rm an}(p)\in\m_p^2\},
 \]
 and ${\rm Zero}(f^{\rm an})$ is the zero locus  of $f^{\rm an}$.

\subsection{Main results}
A data $(X,L,W)$ is called  {\it Landau-Ginzburg model}, or just {\it LG-model},  if  $X$ is a scheme, $L$ is a line bundle  on $X$, and   $W\in\Gamma(X,L)$ is a section of $L$. To a LG-model $(X,L,W)$ we associate a triangulated category $\DMF(X,L,W)$, called the {\it derived matrix factorization category}, introduced by Positselski \cite{posi,efi-posi}. Tensor products of matrix factorizations defines the bifunctor;
\[
(-)\otimes(-):\DMF(X,L,W_1)\times\DMF(X,L,W_2)\to \DMF(X,L,W_1+W_2).
\]
In particular, $\DMF(X,L,W)$  has a tensor action from $\DMF(X,L,0)$.

The following is our main result of classification of thick subcategories of  derived matrix factorization categories.

\begin{thm}[Theorem \ref{main result}]\label{main result in intro}

Let $X$ be a separated Noetherian scheme with an ample family of line bundles, $L$ be a line bundle on $X$, and let $W\in \Gamma(X,L)$ be a non-zero-divisor. Denote by $X_0$ the zero scheme of $W$. Then there is a bijective correspondence
$${\small
\left\{
\begin{aligned}
unions \hspace{2mm}of\hspace{2mm} closed \hspace{4mm}\\
subsets  \hspace{1mm}of \hspace{1mm}{\rm Sing}(X_0/X)\hspace{0mm}
\end{aligned}
\right\}
\begin{aligned}
\xrightarrow{\sigma}\\
\xleftarrow{\tau}\\
\end{aligned}
\left\{
\begin{aligned}
thick \hspace{1.5mm}subcategories\hspace{1.5mm}
 of \hspace{1mm}{\rm DMF}(X,L,W)\hspace{1mm}that \hspace{1mm}are\\
 closed \hspace{1mm}under \hspace{1mm}tensor \hspace{1mm}action \hspace{1mm}from \hspace{1mm}{\rm DMF}(X,L,0)
\end{aligned}
\right\}
}$$
The bijective map $\sigma$ sends $Y$ to the thick subcategory consisting of matrix factorizations $F\in {\rm DMF}(X,L,W)$ with ${\rm Supp}(F)\subseteq Y$.
The inverse bijection $\tau$ sends  $\mathcal{T}$  to the union $\bigcup_{F\in\mathcal{T}}{\rm Supp}(F)$.
\end{thm}

If $X$ is a regular separated Noetherian scheme,  then $X$ has an ample family of line bundles and ${\rm DMF}(X,L,W)$ is equivalent to the singularity category ${\rm D}^{\rm sg}(X_0)$. Furthermore, if $X=\Spec R$ is affine with $R$ regular local, ${\rm DMF}(X,L,W)$ is equivalent to the stable category $\underline{{\rm CM}}(R/W)$ of maximal Cohen-Macaulay modules over the hypersurface $R/W$. Hence Theorem \ref{main result in intro}  can be considered as  a simultaneous generalization  of Stevenson's result in \cite{stev} and Takahashi's result in \cite{takahashi}. 

As an application of the above main result, we see that the closedness of the relative singular locus $\Sing(X_0/X)$ is related to the existence of a $\otimes$-generator of $\DMF(X,L,W)$, where we say that an object $G\in\DMF(X,L,W)$ is a {\it $\otimes$-generator} if the smallest thick subcategory that is closed under tensor action from $\DMF(X,L,0)$ and contains $G$ is $\DMF(X,L,W)$. 

\begin{cor}[Corollary \ref{generator}]
Notation is same as in Theorem \ref{main result in intro}. Then the subset $\Sing(X_0/X)$ of $X_0$ is  closed  if and only if  $\DMF(X,L,W)$ has a $\otimes$-generator.
\end{cor}

Furthermore, we construct the relative singular loci from  the derived matrix factorization categories. If $2\in\Gamma(X,\cO_X)$ is a unit in the ring $\Gamma(X,\cO_X)$, the derived matrix factorization category $\DMF(X,L,W)$ has a natural (pseudo) tensor triangulated structure $\otimes^{\frac{1}{2}}$ on it. Using Theorem \ref{main result in intro} and the theory of Balmer's tensor triangular geometry, 
we prove that  the spectrum of the (pseudo) tensor triangulated category $(\DMF(X,L,W),\otimes^{\frac{1}{2}})$ is the relative singular locus $\Sing(X_0/X)$.

\begin{thm}[Corollary \ref{reconst}]\label{main result 2 in intro}
Let $X$ be a separated Noetherian scheme with an ample ample family of line bundles, and let $W\in \Gamma(X,L)$ be a non-zero-divisor of a line bundle $L$. Assume that $2\in\Gamma(X,\cO_X)$ is a unit. Then we have a homeomorphism
$${\rm Spc}({\rm DMF}(X,L,W), \otimes^{\frac{1}{2}})\cong{\rm Sing}(X_0/X).$$
\end{thm}

This result is a generalization of Yu's result \cite[Theorem 1.2]{yup}, where he proved Theorem \ref{main result 2 in intro} in the case when $X$ is an affine regular scheme of finite Krull dimension  by using the classification result  due to Walker.

\subsection{Plan of the paper}
In section 2 we provide basic definitions  and properties about derived matrix factorization categories.
In section 3 we give the definitions of globally/locally relative singular loci and prove some properties about relative singular loci for zero schemes of regular sections of line bundles. 
In section 4 we prove tensor nilpotence properties of matrix factorizations which are  key properties for our classification result.
In section 5 we prove the main result Theorem \ref{main result in intro}. In section 6 we recall the theory of Balmer's tensor triangular geometry, and we study the natural tensor triangulated structure on derived matrix factorization categories.

\subsection{Acknowledgements} 
The author would like to thank  Hokuto Uehara for many useful comments and his continuous support. The author also thank Michael Wemyss and Shinnosuke Okawa for valuable discussions and comments on a draft version of the paper. The author was a Research Fellow of Japan Society for the Promotion of Science. He was partially supported by Grant-in-Aid for JSPS Fellows No.26-6240.

\vspace{3mm}
\section{Derived matrix factorizations}

\subsection{Derived matrix factorization categories}

In the first subsection, we recall the definition of  the derived matrix factorization category of a Landau-Ginzburg model, which is introduced by Positselski (cf. \cite{posi}, \cite{efi-posi}), and provide its basic properties.

\begin{dfn}
A \textbf{Landau-Ginzburg model}, or \textbf{LG model}, is data $(X,L,W)$ consisting of 
a scheme $X$, an invertible sheaf $L$ on $X$, and  a section $W\in\Gamma(X,L)$ of $L$. 
\end{dfn}

\begin{nota}

If $L$ is isomorphic to the structure sheaf $\mathcal{O}_X$, we denote the LG model by $(X,W)$.  If $X=\Spec R$ is an affine scheme,  we denote the LG model by  $(R,L,W)$, where $L$ is considered as an invertible $R$-module and $W\in L$.
\end{nota}
 
 For a LG model, we consider its factorizations which are $^{\rotatebox[origin=C]{180}{"}}$twisted"  complexes.

\begin{dfn}
Let $(X,L,W)$ be a LG model. A \textbf{factorization} $F$ of $(X,L,W)$ is a sequence
$$F=\Bigl(F_1\xrightarrow{\varphi_1^F} F_0\xrightarrow{\varphi_0^F} F_1\otimes L\Bigr),$$
where each $F_i$ is a coherent sheaf on $X$ and each $\varphi_i^F$ is a  homomorphism such that $\varphi_0^F\circ\varphi_1^F=W\cdot {\rm id}_{F_1}$ and $(\varphi_1^F\otimes L)\circ\varphi_0^F=W\cdot {\rm id}_{F_0}$. Coherent sheaves $F_0$ and $F_1$ in the above sequence are called \textbf{components} of the factorization $F$. If the components $F_i$ of  $F$ are locally free sheaves, we call $F$ a \textbf{matrix factorization} of $(X,L,W)$.
\end{dfn}

\begin{nota}\label{sheaf fact}
We can consider any coherent sheaf $F\in\coh X$ as a factorization of $(X,L,0)$ of the following form
$$\Bigl(0\longrightarrow F\longrightarrow0\Bigr).$$
By abuse of notation, we will often denote the above factorization by the same notation $F$.
\end{nota}

\begin{dfn}
For a LG model $(X,L,W)$, we define an exact category
$${\rm coh}(X,L,W)$$
whose objects are factorizations of $(X,L,W)$, and whose set of morphisms are defined as follows:
For two objects $E,F\in {\rm coh}(X,L,W)$, we define  ${\rm Hom}(E,F)$ as the set of pairs $(f_1,f_0)$ of $f_i\in{\rm Hom}_{{\rm coh}X}(E_i, F_i)$ such that  the following diagram is commutative;
\[\xymatrix{
 E_1\ar[rr]^{\varphi_1^E}\ar[d]_{f_1}&&E_0\ar[rr]^{\varphi_0^E}\ar[d]^{f_0}&&E_1\otimes L\ar[d]^{f_1\otimes L}\\
F_1\ar[rr]^{\varphi_1^F}&&F_0\ar[rr]^{\varphi_0^F}&&F_1\otimes L.
 }\]
 Note that if $X$ is Noetherian, $\coh(X,L,W)$ is an abelian category.
We define a full additive subcategory $${\rm MF}(X,L,W)$$ of  ${\rm coh}(X,L,W)$ whose objects are matrix factorizations. By construction, ${\rm MF}(X,L,W)$ is also an exact category.
\end{dfn}

Since factorizations are $^{\rotatebox[origin=C]{180}{"}}$twisted" complexes, we can consider homotopy category  of factorizations.

\begin{dfn}
Two morphisms $f=(f_1,f_0)$ and $g=(g_1,g_0)$ in ${\rm Hom}_{{\rm coh}(X,L,W)}(E,F)$ are \textbf{homotopy equivalent}, denoted by $f\sim g$, if there exist two homomorphisms in $\coh X$
$$h_0:E_0\rightarrow F_1\hspace{3mm}{\rm and}\hspace{3mm}h_1:E_1\otimes L\rightarrow F_0$$
such that $f_0-g_0=\varphi_1^Fh_0+h_1\varphi_0^E$ and $f_1\otimes L-g_1\otimes L=\varphi_0^Fh_1+(h_0\otimes L)(\varphi_1^E\otimes L)$.
\end{dfn}

\vspace{2mm}
The {\bf homotopy category of factorizations} \[{\rm Kcoh}(X,L,W)\] is defined as the category whose objects are same as ${\rm coh}(X,L,W)$, and the set of morphisms are defined as the set of homotopy equivalence  classes;
$${\rm Hom}_{{\rm Kcoh}(X,L,W)}(E,F):={\rm Hom}_{{\rm coh}(X,L,W)}(E,F)/\sim.$$
Similarly, we define the homotopy category of matrix factorizations ${\rm KMF}(X,L,W)$, i.e. 
$${\rm Ob}({\rm KMF}(X,L,W)):={\rm MF}(X,L,W)$$
$${\rm Hom}_{{\rm KMF}(X,L,W)}(E,F):={\rm Hom}_{{\rm MF}(X,L,W)}(E,F)/\sim.$$

\vspace{3mm}
Next we define the totalization of a bounded complex of factorizations, which is an analogy of the total complex of a double complex.

\begin{dfn}
Let $F^{\text{\tiny{\textbullet}}}=(\cdot\cdot\cdot\rightarrow F^i\xrightarrow{\delta^i}F^{i+1}\rightarrow\cdot\cdot\cdot)$ be a bounded complex of ${\rm coh}(X,L,W)$. For $l=0,1$, set
$$T_l:=\bigoplus_{i+j=-l}F^i_{\overline{j}}\otimes L^{\otimes{\lceil j/2\rceil}},$$
and let $$t_l:T_l\rightarrow T_{\overline{l+1}}$$
be a  homomorphism given by
$$t_l|_{F^i_{\overline{j}}\otimes L^{\otimes{\lceil j/2\rceil}}}:=\delta^i_{\overline{j}}\otimes L^{\otimes{\lceil j/2\rceil}}+(-1)^i\varphi^{F^i}_{\overline{j}}\otimes L^{\otimes{\lceil j/2\rceil}},$$
where $\overline{n}$ is $n$ modulo $2$, and $\lceil m\rceil$ is the minimum integer which is greater than or equal to a real number $m$.
We define the {\bf totalization} Tot$(F^{\text{\tiny{\textbullet}}})\in {\rm coh}(X,L,W))$ of $F^{\text{\tiny{\textbullet}}}$ as 
$${\rm Tot}(F^{\text{\tiny{\textbullet}}}):=\Bigl(T_1\xrightarrow{t_1}T_0\xrightarrow{t_0}T_1\otimes L\Bigr).$$
\end{dfn}

\vspace{3mm}In what follows, we will recall that the homotopy categories ${\rm Kcoh}(X,L,W)$ and ${\rm KMF}(X,L,W)$ have  structures of  triangulated categories.

\vspace*{2mm}
\begin{dfn}
We define an automorphism $T$ on ${\rm Kcoh}(X,L,W))$, which is called \textbf{shift functor}, as follows.
For an object $F\in {\rm Kcoh}(X,L,W)$, we define an object $T(F)$ as
$$T(F):=\Bigl(F_0\xrightarrow{-\varphi^F_0}F_1\otimes L\xrightarrow{-\varphi^F_1\otimes L}F_0\otimes L\Bigr),$$
and for a morphism $f=(f_1,f_0)\in {\rm Hom}(E,F)$ we set $T(f):=(f_0,f_1\otimes L)\in{\rm Hom}(T(E),T(F))$. For any integer $n\in\mathbb{Z}$, denote by $(-)[n]$ the functor $T^n(-)$.
\end{dfn}

\vspace*{2mm}
\begin{dfn}
Let $f : E\rightarrow F$ be a morphism in ${\rm coh}(X,L,W)$. We define its {\bf mapping cone} Cone$(f)$ to be the totalization of the complex $$(\cdot\cdot\cdot\rightarrow0\rightarrow E\xrightarrow{f} F\rightarrow0\rightarrow\cdot\cdot\cdot)$$ with $F$ in degree zero.

A $\textbf{distinguished triangle}$ is a sequence in ${\rm Kcoh}(X,L,W)$ which is isomorphic to a sequence of the form
$$E\xrightarrow{f}F\xrightarrow{i}{\rm Cone}(f)\xrightarrow{p}E[1],$$
where $i$ and $p$ are  natural injection and  projection respectively.
\end{dfn}

The following proposition is well known to experts.

\begin{prop}
The homotopy categories ${\rm Kcoh}(X,L,W)$ and ${\rm KMF}(X,L,W)$ are triangulated categories with respect to the above shift functor and the above distinguished triangles.
\end{prop}

\vspace*{2mm}
Following Positselski (\cite{posi}, \cite{efi-posi}), we define  derived factorization categories.

\begin{dfn}
Denote by  ${\rm Acoh}(X,L,W)$ the smallest thick subcategory of ${\rm Kcoh}(X,L,W)$ containing all totalizations of short exact sequences in ${\rm coh}(X,L,W)$.  We define the \textbf{derived factorization category} of $(X,L,W)$ as the Verdier quotient 
$${\rm Dcoh}(X,L,W):={\rm Kcoh}(X,L,W)/{\rm Acoh}(X,L,W).$$
Similarly, we consider the thick subcategory ${\rm AMF}(X,L,W)$ containing all totalizations of short exact sequences in the exact category ${\rm MF}(X,L,W)$, and define the \textbf{derived matrix factorization category} by
$${\rm DMF}(X,L,W):={\rm KMF}(X,L,W)/{\rm AMF}(X,L,W).$$

\end{dfn}

\vspace{2mm}

The following proposition is a special case of \cite[Lemma 2.24]{bdfik}.

\begin{prop}[cf. {\cite[Lemma 2.24]{bdfik}}]\label{affine}
Assume that $X=\Spec R$ is an affine scheme. For $P\in {\rm KMF}(R,L,W)$ and $A\in{\rm Acoh}(R,L,W)$, we have 
$${\rm Hom}_{{\rm Kcoh}(R,L,W)}(P,A)=0.$$
In particular, the Verdier  localizing functor 
$${\rm KMF}(R,L,W)\xrightarrow{\sim} {\rm DMF}(R,L,W)$$
is an equivalence.

\end{prop}

For later use, we consider larger categories of factorizations. Denote by ${\rm Sh}(X,L,W)$ the abelian category whose objects are factorizations whose components are $\mathcal{O}_X$-modules. More precisely, objects of  ${\rm Sh}(X,L,W)$ are sequences of the following form
$$F=\Bigl(F_1\xrightarrow{\varphi_1^F} F_0\xrightarrow{\varphi_0^F} F_1\otimes L\Bigr),$$
where $F_i$ are $\mathcal{O}_X$-modules  and $\varphi_i^F$ are  homomorphisms such that $\varphi_0^F\circ\varphi_1^F=W\cdot {\rm id}_{F_1}$ and $\varphi_1^F\otimes L\circ\varphi_0^F=W\cdot {\rm id}_{F_0}$. 
Denote by ${\rm Qcoh}(X,L,W)$, ${\rm InjSh}(X,L,W)$, and ${\rm InjQcoh}(X,L,W)$ the full subcategories of ${\rm Sh}(X,L,W)$ consisting of factorizations whose components are quasi-coherent sheaves, injective $\mathcal{O}_X$-modules, and injective quasi-coherent sheaves respectively.

Then, similarly to  ${\rm Kcoh}(X,L,W)$, we can consider their homotopy categories ${\rm KSh}(X,L,W)$, ${\rm KQcoh}(X,L,W)$,  ${\rm KInjSh}(X,L,W)$,  ${\rm KInjQcoh}(X,L,W)$ respectively,  and these homotopy categories have  natural triangulated structures similar to ${\rm Kcoh}(X,L,W)$.

\begin{dfn}
Denote by ${\rm A^{co}Sh}(X,L,W)$ (resp. ${\rm A^{co}Qcoh}(X,L,W)$) the smallest thick subcategory of  ${\rm KSh}(X,L,W)$ (resp. ${\rm KQcoh}(X,L,W)$) containing all totalizations of short exact sequences in ${\rm Sh}(X,L,W)$ (resp. ${\rm Qcoh}(X,L,W)$) and closed under arbitrary direct sums. Following \cite{posi}, \cite{efi-posi}, we define  the \textbf{coderived factorization categories} ${\rm D^{co}Sh}(X,L,W)$ and ${\rm D^{co}Qcoh}(X,L,W)$ as the following Verdier quotients 
\[{\rm D^{co}Sh}(X,L,W):={\rm KSh}(X,L,W)/{\rm A^{co}Sh}(X,L,W)\]
\[{\rm D^{co}Qcoh}(X,L,W):={\rm KQcoh}(X,L,W)/{\rm A^{co}Qcoh}(X,L,W).\]

\end{dfn}

\begin{lem}[{\cite{bdfik}}, \cite{efi-posi}]\label{ff} Assume that $X$ is Noetherian.
\begin{itemize}
\item[$(1)$]The natural functor ${\rm KInjSh}(X,L,W)\rightarrow{\rm D^{co}Sh}(X,L,W)$ is  an equivalence.
\item[$(2)$] The natural functor ${\rm KInjQcoh}(X,L,W)\rightarrow{\rm D^{co}Qcoh}(X,L,W)$ is  an equivalence.
\item[$(3)$] The natural functor ${\rm D^{co}Qcoh}(X,L,W)\rightarrow{\rm D^{co}Sh}(X,L,W)$ is fully faithful.
\item[$(4)$] The natural functor ${\rm Dcoh}(X,L,W)\rightarrow{\rm D^{co}Qcoh}(X,L,W)$ is fully faithful.
\item[$(5)$] The natural functor ${\rm DMF}(X,L,W)\rightarrow{\rm Dcoh}(X,L,W)$ is fully faithful.

\end{itemize}

\begin{proof}
$(1)$ and $(2)$ follow from \cite[Corollary 2.25]{bdfik}. $(3)$ follows from $(1)$ and $(2)$. $(4)$ and $(5)$ are \cite[Propostion 1.5.(d)]{efi-posi} and \cite[Corollary 2.3.(i)]{efi-posi} respectively.
\end{proof}
\end{lem}

\subsection{Case when $W=0$}
In this section, we consider cases when $W=0$. Firstly, we will define cohomologies of  factorizations of $(X,L,0)$. 

\begin{dfn}
For an object $F\in{\rm Qcoh}(X,L,0)$, we define its cohomologies ${\sf H}_i(F)\in \Qcoh X$ as 
$${\sf H}_i(F):={\rm Ker}(\varphi_i^F)/{\rm Im}(\varphi_{i-1}^F)\hspace{6mm} {\rm for }\hspace{3mm}{i\in \mathbb{Z}/2\mathbb{Z}}$$

\end{dfn}

\begin{lem}\label{field split}
Let $k$ be any field. Then, for any object $F\in{\rm KMF}(k,0)$, there are two finite dimensional $k$-vector spaces $V_1$ and $V_2$ such that $F$ is isomorphic to $V_1\oplus V_2[1]$ in ${\rm KMF}(k,0)$, where $V_i$ denotes the factorization of the form $(0\rightarrow V_i\rightarrow0)$ by Notation \ref{sheaf fact}.
\begin{proof}
By \cite[Lemma 2.26]{bdfik}, there are  two finite dimensional $k$-vactor spaces $V$ and $V'$, and a triangle   of the following form in ${\rm Dcoh}(k,0)={\rm DMF}(k,0)$
$$V\rightarrow V'\rightarrow F\rightarrow V[1].$$
But ${\rm DMF}(k,0)={\rm KMF}(k,0)$ by Proposition \ref{affine}, so we have a $k$-linear homomorphism $f:V\rightarrow V'$ such that $F$ is isomorphic to $C(\overline{f})$, where $\overline{f}$ is the morphism in ${\rm KMF}(k,0)$ represented by the following morphism in ${\rm MF}(k,0)$
\[\xymatrix{
0\ar[rr]\ar[d]&&V\ar[d]^{f}\ar[rr]&&\ar[d]0\\
0\ar[rr]&&V'\ar[rr]&&0
}\]
By construction of mapping cones, $C(\overline{f})$ is isomorphic to the following matrix factorization
$$\Bigl( V\xrightarrow{f}V'\xrightarrow{0}V\Bigr).$$ 
Let $I:={\rm Im}(f)$ be the image of $f$, and let $K:={\rm Ker}(f)$ be the kernel of $f$. Then there is a $k$-vector space $J$ such that $V'=I\oplus J$. Since $V=K\oplus I$, we have the following isomorphism in ${\rm MF}(k,0)$ 
$$F\cong \Bigl( K\xrightarrow{0}0\xrightarrow{0}K\Bigr)\oplus\Bigl( I\xrightarrow{\sim}I\xrightarrow{0}I\Bigr)\oplus\Bigl( 0\xrightarrow{0}J\xrightarrow{0}0\Bigr)$$
But the object $\Bigl( I\xrightarrow{\sim}I\xrightarrow{0}I\Bigr)$ is zero in ${\rm KMF}(k,0)$. Hence $F\cong J\oplus K[1]$ in ${\rm KMF}(k,0)$.
\end{proof}
\end{lem}

\begin{cor}\label{stalk lem} 
Let $k$ be a field. Any non-zero morphism $f:(0\to k\to 0)\rightarrow E$  in ${\rm KMF}(k,0)$ is  a split mono.

\begin{proof}
This follows   from Lemma \ref{field split}.
\end{proof}
\end{cor}

\subsection{Tensor products and sheaf Homs functors}

In this subsection, we recall tensor products and local homs  on derived matrix factorization categories. Let $(X,L,W)$ be a LG model, and  $V\in\Gamma(X,L)$ be another global section.

For $E\in{\rm MF}(X,L,V)$ and $F\in{\rm MF}(X,L,W)$, we define the tensor product $$E\otimes F\in{\rm MF}(X,L,V+W)$$ of $E$ and $F$ as 
$$(E\otimes F)_1:=\bigl(E_1\otimes F_0\bigr)\oplus \bigl(E_0\otimes F_1\bigr),$$
$$(E\otimes F)_0:=\bigl(E_0\otimes F_0\bigr)\oplus \bigl(E_1\otimes F_1\otimes L\bigr),$$
\,$${\varphi_1^{E\otimes F}}:=\begin{pmatrix}
\varphi^E_1\otimes1&1\otimes\varphi^F_1\\ -1\otimes\varphi^F_0&\varphi^E_0\otimes1
\end{pmatrix},
$$
and 
$$\varphi_0^{E\otimes F}:=\begin{pmatrix}
\varphi^E_0\otimes1 & -1\otimes\varphi^F_1\\
1\otimes\varphi^F_0 & \varphi^E_1\otimes1
\vspace{2mm}
\end{pmatrix}.$$
This defines an additive functor $(-)\otimes(-):{\rm MF}(X,L,V)\times{\rm MF}(X,L,W)\rightarrow{\rm MF}(X,L,V+W)$, and it naturally induces an exact functor $$(-)\otimes(-):{\rm DMF}(X,L,V)\times{\rm DMF}(X,L,W)\rightarrow{\rm DMF}(X,L,V+W).$$

\vspace{4mm}
We define the sheaf Hom 
$$\mathcal{H}om(E,F)\in{\rm MF}(X,L,W-V)$$
from $E$ to $F$ as 
$$\mathcal{H}om(E,F)_1:=(\mathcal{H}om(E_1,F_0)\otimes L^{-1})\oplus\mathcal{H}om(E_0,F_1),$$
$$\mathcal{H}om(E,F)_0:=\mathcal{H}om(E_0,F_0)\oplus\mathcal{H}om(E_1,F_1),$$
\,$${\varphi_1^{\mathcal{H}om(E,F)}}:=\begin{pmatrix}
(\ast)\circ\varphi^E_0&\varphi^F_1\circ(*)\\ 
(\varphi^F_0\otimes L^{-1})\circ(*)&(*)\circ\varphi^E_1
\end{pmatrix},
$$
and 
$$\varphi_0^{\mathcal{H}om(E,F)}:=\begin{pmatrix}
-(*)\circ\varphi^E_1 &\varphi^F_1\circ(*)\\
\varphi^F_0\circ(*) & -(*)\circ(\varphi^E_0\otimes L^{-1})
\end{pmatrix}.\vspace{2mm}$$
This defines an additive functor $\mathcal{H}om(-,-):{\rm MF}(X,L,V)^{\rm op}\times{\rm MF}(X,L,W)\rightarrow{\rm MF}(X,L,W-V)$, and it induces an exact functor 
$$\mathcal{H}om(-,-):{\rm DMF}(X,L,V)^{\rm op}\times{\rm DMF}(X,L,W)\rightarrow{\rm DMF}(X,L,W-V).$$

\vspace{3mm}
The following is  standard, so we skip the proof (see \cite{bfk} or \cite{ls} for details) . 
\begin{prop}\label{hom and tens}
Let $E\in{\rm MF}(X,L,V)$, $F\in{\rm MF}(X,L,W)$, and $G\in{\rm MF}(X,L,V+W)$.
\begin{itemize}
\item[$(1)$]
 We have a natural isomorphism 
$${\rm Hom}_{{\rm DMF}(X,L,V+W)}(E\otimes F,G)\cong{\rm Hom}_{{\rm DMF}(X,L,V)}(E,\mathcal{H}om(F,G)).$$ 
\item[$(2)$]
There is a natural isomorphism in ${\rm MF}(X,L,V+W)$
$$\mathcal{H}om(G,E\otimes F)\cong \mathcal{H}om(G,E)\otimes F.$$
\end{itemize}
\end{prop}

Recall that $\mathcal{O}_X\in{\rm MF}(X,L,0)$ denotes the matrix factorization of the form  $\Bigl(0\rightarrow\mathcal{O}_X\rightarrow0\Bigr)$ by Notation \ref{sheaf fact}. For any object $F\in{\rm MF}(X,L,W)$, we define \textbf{the dual} $$F^{\vee}:=\mathcal{H}om(F,\mathcal{O}_X)\in{\rm MF}(X,L,-W)$$ of $F$. By Proposition \ref{hom and tens}, the functors
$(-)\otimes F:{\rm DMF}(X,L,V)\rightarrow {\rm DMF}(X,L,V+W)$ and $(-)\otimes F^{\vee}:{\rm DMF}(X,L,V+W)\rightarrow {\rm DMF}(X,L,V)$ are adjoint;
$$(-)\otimes F\dashv(-)\otimes F^{\vee}.$$

\subsection{Supports of matrix factorizations}

We study the supports of objects in derived matrix factorization categories. Let $(X,L,W)$ be a LG model.
  
  For any point $p\in X$, we denote  by $X_p:=\Spec(\mathcal{O}_{X,p})$ the stalk of $X$ at $p$, and let ${\rm vect}\,X$ be the category of locally free sheaves of finite ranks on $X$.  
Taking the stalk $(-)_p: {\rm vect}\,X\to {\rm vect}\, X_p$ at $p$ induces the functor $(-)_p:\KMF(X,L,W) \to \KMF(X_p,W_p)$ defined by $$F_p:=\Bigl((F_1)_p\xrightarrow{(\varphi_1)_p}(F_0)_p\xrightarrow{(\varphi_0)_p}(F_1)_p\Bigr).$$ 
Since the functor $(-)_p:{\rm vect}\, X\to {\rm vect}\, X_p$ preserves short exact sequences, the induced functor  $(-)_p:\KMF(X,L,W) \to \KMF(X_p,W_p)$ maps $\AMF(X,L,W)$ to $\AMF(X_p,W_p)$. Hence it defines the following functor;
$$(-)_p:{\rm DMF}(X,L,W)\rightarrow {\rm KMF}(X_p, W_p),$$
since we have the natural equivalence ${\rm KMF}(X_p, W_p)\xrightarrow{\sim}{\rm DMF}(X_p, W_p)$ by Proposition \ref{affine}.

\begin{dfn}\label{def supp}
For an object $F\in {\rm DMF}(X,L,W)$, we define its  support as
$${\rm Supp}(F):=\{\,p\in X\,\mid\, F_p\ne0 \in {\rm KMF}(X_p, W_p)\}.$$
\end{dfn}

\begin{prop}\label{supp local} Let $F\in {\rm DMF}(X,L,W)$ be an object.
 \begin{itemize}
 
 \item[$(1)$] If $X=\bigcup_{i\in I}U_i$ is an open covering of $X$, we have the equality of subsets of $X$
$${\rm Supp}(F)=\bigcup_{i\in I}{\rm Supp}(F|_{U_i}),$$
where $F|_{U_i}$ is the restriction of $F$ to ${\rm DMF}(U_i,L|_{U_i},W|_{U_i})$. 

\item[$(2)$] ${\rm Supp}(F)$ is a closed subset of $X$.

\end{itemize}
\begin{proof}
(1)  This follows from isomorphisms $F_p\cong(F|_{U_i})_p$ for any $p\in U_i$.\\
(2) We show the following equality $${\rm Supp}(F)^{\rm c}=\bigcup_{U\in\mathcal{U}}U,$$
where $\mathcal{U}:=\{$ $U \mid U$ is an open subscheme of $X$ such that $F|_U=0$ in ${\rm DMF}(U,L|_U,W|_U)\}$.  The inclusion ${\rm Supp}(F)^{\rm c}\supset\bigcup_{U\in\mathcal{U}}U$ is obvious. We verify that ${\rm Supp}(F)^{\rm c}\subseteq\bigcup_{U\in\mathcal{U}}U$. By definition, for any $p\in {\rm Supp}(F)^{\rm c}$, $F_p=0$ in  ${\rm KMF}(X_p, W_p)$. Let $h=(h_1,h_0)$ be a homotopy giving the homotopy equivalence ${\rm id}_{F_p}\sim 0$. Then there is a neighborhood $U$ of $p$ such that there exist morphisms  $\overline{h_1}:F_1|_U\otimes L_U\rightarrow F_0|_U$ and $\overline{h_0}:F_0|_U\rightarrow F_1|_U$ in $\coh U$ with $\overline{h_1}_p=h_1$ and $\overline{h_0}_p=h_0$. Furthermore, since ${\rm id}_{(F_0)_p}-(\varphi_1)_ph_0-h_1(\varphi_0)_p=0$ in $\coh U_p$, there exists  an open neighborhood  $V\subseteq U$ of $p$ such that ${\rm id}_{F_0|_V}-\varphi_1|_V\overline{h_1}|_V-\varphi_0|_V\overline{h_0}|_V=0$. Then $h_V=(\overline{h_1}|_V,\overline{h_0}|_V)$ gives a homotopy equivalence ${\rm id}_{F|_V}\sim0$. Hence $F|_V=0$ in ${\rm KMF}(V,L|_V,W|_V)$, in particular, so is in ${\rm DMF}(V,L|_V,W|_V)$. Therefore, we have $V\in \mathcal{U}$, which implies that $p\in \bigcup_{U\in\mathcal{U}}U$.
\end{proof}
\end{prop}

\begin{dfn}
Let $F\in {\rm DMF}(X,L,W)$ be an object. For any point $p\in X$, let $\iota_p:\Spec k(p)\rightarrow X$ be a natural  morphism, where $k(p):=\mathcal{O}_{X,p}/\mathfrak{m}_p$ is the residue field of the local ring $(\mathcal{O}_{X,p},\m_p)$. Then we denote by $W\otimes k(p)$ the pull-back $\iota_p^*W$, and we set 
$$F\otimes k(p):=\iota_p^*F=\Bigl(   \iota_p^*F_1\xrightarrow{\iota_p^*\varphi_1} \iota_p^*F_0\xrightarrow{\iota_p^*\varphi_0}\iota_p^*F_1   \Bigr)\in {\rm KMF}(k(p),W\otimes k(p)).$$ 
Then $(-)\otimes k(p)$ defines an exact functor 
$$(-)\otimes k(p):{\rm DMF}(X,L,W)\rightarrow  {\rm KMF}(k(p),W\otimes k(p)).$$

\end{dfn}

The following lemma is a version of Nakayama's lemma for matrix factorizations.

\begin{lem}\label{stalk lem2} Let $F\in {\rm DMF}(X,L,W)$ be an object, and let $p\in X$ be a point. 
 Then $$p\in {\rm Supp}(F) \hspace{3mm}{\rm if  \hspace{1mm}and \hspace{1mm} only \hspace{1mm} if} \hspace{3mm} F\otimes k(p)\ne0  \hspace{2mm}{\rm in}  \hspace{2mm}{\rm KMF}(k(p),W\otimes k(p)).$$

\begin{proof}
Recall that we have natural equivalences ${\rm KMF}(k(p),W\otimes k(p))\xrightarrow{\sim} {\rm DMF}(k(p),W\otimes k(p))$ and ${\rm KMF}(X_p,W_p)\xrightarrow{\sim}{\rm DMF}(X_p,W_p)$. Note that the following diagram of functors is commutative;
\[\xymatrix{
 {\rm DMF}(X,L,W)\ar[rrrr]^{(-)\otimes k(p)}\ar[rrd]_{(-)_p}&&&&{\rm KMF}(k(p),W\otimes k(p))\\
&&{\rm KMF}(X_p,W_p)\ar[rru]_{(-)\otimes k(\mathfrak{m}_p)}&&
 }\]
 where $\mathfrak{m}_p\in X_p$ is the unique closed point.
Hence, if $F_p=0$ in ${\rm KMF}(X_p,W_p)$, then $F\otimes k(p)=0$ in ${\rm KMF}(k(p),W\otimes k(p))$.

For the other implication, it suffices to show that for a local ring $(R,\mathfrak{m})$, an element $w\in R$, and an object $E\in {\rm KMF}(R,w)$, if $E\otimes_R R/\mathfrak{m}=0$ in ${\rm KMF}(R/\mathfrak{m},w\otimes R/\mathfrak{m})$, then $E=0$ in ${\rm KMF}(R,w)$. Since $R$ is local, any locally free modules are free. Hence, the object $E$ can be represented by some matrix factorization of the following form
$$\Bigl(R^{\oplus n_1}\xrightarrow{\varphi_1} R^{\oplus n_0}\xrightarrow{\varphi_0} R^{\oplus n_1} \Bigr).$$
If $E\otimes_R R/\mathfrak{m}=0$, there exist homotopies $h_0:(R/\mathfrak{m})^{\oplus n_0}\rightarrow (R/\mathfrak{m})^{\oplus n_1}$ and $h_1:(R/\mathfrak{m})^{\oplus n_1}\rightarrow (R/\mathfrak{m})^{\oplus n_0}$ such that ${\rm id}_{(R/\mathfrak{m})^{\oplus n_0}}=(\varphi_1\otimes R/\mathfrak{m})h_0+h_1(\varphi_0\otimes R/\mathfrak{m})$ and ${\rm id}_{(R/\mathfrak{m})^{\oplus n_1}}=(\varphi_0\otimes R/\mathfrak{m})h_1+h_0(\varphi_1\otimes R/\mathfrak{m})$. Since $h_0$ and $h_1$ can be represented by  a matrix of units in $R$, there exist homomorphisms $\overline{h_0}:R^{\oplus n_0}\rightarrow R^{\oplus n_1}$ and $\overline{h_1}:R^{\oplus n_1}\rightarrow R^{\oplus n_0}$ such that $\overline{h_i}\otimes_R R/\mathfrak{m}=h_i$ for $i=0,1$. Set
$$\alpha_1:=\varphi_0\overline{h_1}+\overline{h_0}\varphi_1:R^{\oplus n_1}\longrightarrow R^{\oplus n_1}$$
$$\alpha_0:=\varphi_1\overline{h_0}+\overline{h_1}\varphi_0:R^{\oplus n_0}\longrightarrow R^{\oplus n_0}.$$
Then the pair $\alpha:=(\alpha_1,\alpha_0)$ defines an endomorphism of $E$ in the exact category ${\rm MF}(R,w)$. By construction, $\alpha=0$ in ${\rm KMF}(R,w)$. To show that $E=0$ in ${\rm KMF}(R,w)$, it is enough to show that $\alpha:E\rightarrow E$ is an automorphism in ${\rm MF}(R,w)$. For each $i\in \{0,1\}$, we only need to show that $\alpha_i$ is an automorphism. Since  the tensor product $(-)\otimes_R R/\mathfrak{m}$ is a right exact functor and $\alpha_i\otimes_R R/\mathfrak{m}={\rm id}$, we have $${\rm Cok}(\alpha_i)\otimes_R R/\mathfrak{m}\cong{\rm Cok}(\alpha_i\otimes_R R/\mathfrak{m})=0.$$  
By Nakayama's lemma, the above implies that ${\rm Cok}(\alpha_i)=0$. Hence $\alpha_i$ is an automorphism by \cite[Theorem 2.4]{matsumura}. 
\end{proof}
\end{lem}

For later use, we provide the following lemma.

\begin{lem}\label{koszul supp}
Let $R$ be a ring, and let $F=\Bigl(R\xrightarrow{f_1}R\xrightarrow{f_0}R\Bigr)\in{\rm KMF}(R,f_0f_1)$ be an object.
Then we have the following equality of subsets of $\Spec R$:
$${\rm Supp}(F)=\bigcap_{i=0,1}Z(f_i),$$
where $Z(f_i):=\{\,p\in\Spec R\,|\,f_i\otimes k(p)=0\,\}$.
\begin{proof}
($\subseteq$) If $p\in{\rm Supp}(F)$, then $F\otimes k(p)\neq0$ in ${\rm KMF}(k(p), (f_0f_1)\otimes k(p))$ by Lemma \ref{stalk lem2}.
Suppose $p\notin Z(f_i)$ for some $i$. Then  $f_i\otimes k(p)$ is a unit in $k(p)$, and hence  $F\otimes k(p)$ is homotopic to zero. This contradicts to $F\otimes k(p)\neq0$. Hence $p\in \bigcap_{i=0,1}Z(f_i)$.\\
($\supseteq$) If $p\in \bigcap_{i=0,1}Z(f_i)$. Then $f_i\otimes k(p)=0$ for $i=0,1$, and so ${\sf H}_i(F\otimes k(p))=k(p)\neq0$. Hence by \cite[Proposition 2.30]{ls}, $F\otimes k(p)\neq0$. Again by Lemma \ref{stalk lem2}, $p\in {\rm Supp}(F)$.
\end{proof}
\end{lem}

The following lemma is useful to compute the support of tensor products of matrix factorizations.

\begin{lem}\label{supp tens}
Let $V,W\in\Gamma(X,L)$ be any global sections of $L$, and let $E\in{\rm DMF}(X,L,V)$ and $F\in{\rm DMF}(X,L,W)$ be objects. Let $p\in X$ be a point such that $V\otimes k(p)=W\otimes k(p)=0$. Then $p\in {\rm Supp}(E\otimes F)$ if and only if $p\in{\rm Supp}(E)\cap{\rm Supp}(F)$.
\begin{proof}
We have $(E\otimes F)\otimes k(p)\cong (E\otimes k(p))\otimes (F\otimes k(p))$ in ${\rm KMF}(k(p),0)$. By Lemma \ref{stalk lem2}, it is enough to show that  
for a field $k$, and for objects $M,N\in{\rm KMF}(k,0)$, $M\otimes N\neq0$ if and only if $M\neq0$ and $N\neq0$. By Lemma \ref{field split}, for $i=0,1$, we may assume $\varphi_i^M=\varphi_i^N=0$, and then ${\sf H}_i(M)=M_i$ and ${\sf H}_i(N)=N_i$. Then, since $\varphi_i^{M\otimes N}=0$ for $i=0,1$, we have 
\begin{center}${\sf H}_1(M\otimes N)=\Bigl({\sf H}_1(M)\otimes {\sf H}_0(N)\Bigr)\oplus \Bigl({\sf H}_0(M)\otimes {\sf H}_1(N)\Bigr)$\end{center}\!\!\!\!\!\!\!
\begin{center}${\sf H}_0(M\otimes N)=\Bigl({\sf H}_0(M)\otimes {\sf H}_0(N)\Bigr)\oplus \Bigl({\sf H}_1(M)\otimes {\sf H}_1(N)\Bigr).$\end{center}
Hence, by \cite[Proposition 2.30]{ls}, we see that $M\otimes N\neq0$ if and only if $M\neq0$ and $N\neq0$.
\end{proof}
\end{lem}

At the end of this section, we organize  fundamental  properties of supports of matrix factorizations.

\begin{lem}\label{supp data} Let $E,F,G\in{\rm DMF}(X,L,W)$ be objects. We have the following.

$(1)$ ${\rm Supp}(E\oplus F)={\rm Supp}(E)\cup{\rm Supp}(F)$.

$(2)$ ${\rm Supp}(F[1])={\rm Supp}(F)$.

$(3)$ ${\rm Supp}(E)\subseteq{\rm Supp}(F)\cup{\rm Supp}(G)$ for any distinguished triangle $E\rightarrow F\rightarrow G\rightarrow E[1]$.

$(4)$ ${\rm Supp}(E\otimes F)={\rm Supp}(E)\cap{\rm Supp}(F)$.
\begin{proof}
(1), (2), and (3) are obvious.   If $p\in {\rm Supp}(M)$ for some object $M\in {\rm DMF}(X,L,W)$, then $W\otimes k(p)=0$, since ${\rm KMF}(k(p),W\otimes k(p))=0$ if $W\otimes k(p)\neq0$. Hence (4) follows from Lemma \ref{supp tens}.
\end{proof}
\end{lem}

\section{Relative singular locus and singularity category} 
In this section,  we define  relative singular loci and prove some properties about it.
Let $S$ be a Noetherian scheme and let $F\in {\rm D^b}(\coh S)$ be a bounded complex of coherent sheaves. The complex $F$ is called \textbf{perfect} if it is locally quasi-isomorphic to a bounded complex of locally free sheaves of finite rank. $\Perf S\subset{\rm D^b}(\coh S)$ denotes the thick subcategory of perfect complexes.

 We define globally/locally relative singular locus. Recall our notation $S_p:=\Spec({\mathcal{O}_{S,p}})$ for 
 any point $p\in S$.
 
\begin{dfn}
Let $S$ be a Noetherian scheme,  and let $i:T\hookrightarrow S$ be a closed immersion.  

\begin{itemize}
\item[$(1)$] The subset ${\rm Sing}(T/S)\subset T$,  called  the \textbf{singular locus  of $T$ globally relative} (or just \textbf{relative}) \textbf{to  $S$}, is defined by 
$$\Sing(T/S):=\{\,p\in T\,\mid\, \exists F\in \coh T \hspace{2mm}{\rm such\hspace{1mm} that}\hspace{2mm} F_p\notin\Perf T_p \hspace{2mm}{\rm and}\hspace{2mm} i_*(F)\in \Perf S\,\}$$

\item[$(2)$] The subset  ${\rm Sing^{loc}}(T/S)\subset T$,  called the \textbf{singular locus   of $T$ locally relative to $S$}, is defined by 
\[\Singloc(T/S):=\{\,p\in T\,\mid\, \exists F\in \coh T_p \hspace{2mm}{\rm such\hspace{1mm} that}\hspace{2mm} F\notin\Perf T_p \hspace{2mm}{\rm and}\hspace{2mm} {i_p}_*(F)\in \Perf S_p\,\}\]
where $i_p:T_p\hookrightarrow S_p$ is the closed immersion induced by $i:T\hookrightarrow S$
\end{itemize}
\end{dfn}

\vspace{1mm}

\begin{prop}\label{rel sing = sing}

Let $S$ be a Noetherian scheme, and let $i:T\hookrightarrow S$ be a closed immersion. Then we have $${\rm Sing}(T/S)\subseteq {\rm Sing^{loc}}(T/S)\subseteq {\rm Sing}(T).$$
Furthermore, if $S$  is  regular,  globally and locally relative singular loci coincide with usual singular locus; $${\rm Sing}(T/S)={\rm Sing^{loc}}(T/S)={\rm Sing}(T).$$
\begin{proof} 

The first assertion ${\rm Sing}(T/S)\subseteq {\rm Sing^{loc}}(T/S)\subseteq {\rm Sing}(T)$ is obvious.   

For the latter assertion, assume that $S$ is regular. If ${\rm Sing}(T)=\emptyset$, ${\rm Sing}(T/S)={\rm Sing^{loc}}(T/S)={\rm Sing}(T)=\emptyset$ by the former assertion. Assume that ${\rm Sing}(T)\neq\emptyset$, and let $p\in{\rm Sing}(T)$ be a singular  point. It is enough to show that $p\in{\rm Sing}(T/S)$. Since  the projective dimension, denoted by ${\rm pd}_{\mathcal{O}_{T,p}}k(p)$, of $k(p)$ as $\mathcal{O}_{T,p}$-module coincides with the global dimension of $\mathcal{O}_{T,p}$, we have ${\rm pd}_{\mathcal{O}_{T,p}}(k(p))=\infty$. This implies that  $k(p)\notin\Perf T_p$. Let $\{a_1,\hdots,a_r\}\subset\mathcal{O}_{T,p}$ be  generators of the maximal ideal $\mathfrak{m}_p$ of the local ring $\mathcal{O}_{T,p}$. Then there exist a small open affine neighborhood $U=\Spec R\subset T$ of $p$ and elements $b_1,\hdots, b_r\in R$ such that $({b_i})_p=a_i$. Let $I:=\langle b_1,\hdots,b_r\rangle$ be the ideal of $R$ generated by $b_i$. Then $I_p\cong \mathfrak{m}_p$ and  $(R/I)_p\cong k(p)$. Take an extension $F\in \coh T$ of the coherent sheaf $\widetilde{R/I}\in\coh U$, i.e. $F|_U\cong \widetilde{R/I}$. Then $F_p\cong k(p)\notin\Perf T_p$.  Since $S$ is regular, we have  ${\rm D^{b}}(\coh S)=\Perf S$, and so  $i_*(F)\in\Perf S$. Hence $p\in {\rm Sing}(T/S)$.
\end{proof}
\end{prop}

The locally relative singular locus has a local property.
\begin{lem} Let $i:T\hookrightarrow S$ be a closed immersion of Noetherian schemes. Then we have 
\[
\Singloc(T/S)=\bigcup_{p\in T}\Singloc(T_p/S_p),
\]
where the sets $\Singloc(T_p/S_p)$ on the right hand side  are considered as the subsets of $T$ via the natural injective maps $j_p:T_p\hookrightarrow T$.
\begin{proof}
If  $p\in \Singloc(T/S)$, then $\m_p\in\Singloc(T_p/S_p)$, and $j_p(\m_p)=p$. This means that  $p\in\Singloc(T_p/S_p)$, and so $\Singloc(T/S)\subseteq\bigcup_{p\in T}\Singloc(T_p/S_p)$.  If  $q\in\Singloc(T_p/S_p)$  for some $p\in T$, then, since $(T_p)_q\cong T_{j_p(q)}$, $j_p(q)\in\Singloc(T/S)$. Hence $\bigcup_{p\in T}\Singloc(T_p/S_p)\subseteq\Singloc(T/S)$.
\end{proof}
\end{lem}

\vspace{4mm}

Next we recall singularity categories. Let $X$ be a separated Noetherian scheme with {\it resolution property}, i.e.~for any $F\in \coh X$, there exist a  locally free coherent sheaf $E$ and a surjective homomorphism $E\twoheadrightarrow F$. Following \cite{orlov2}, we define the triangulated category of singularities ${\rm D^{sg}}(X)$ as the Verdier quotient \[\Dsg(X)={\rm D^b}(\coh X)/\Perf X.\] In our assumption, $\Perf X$ coincides with thick subcategory of complexes which are quasi-isomorphic to a bounded complex of locally free sheaves of finite rank. 

We recall that derived matrix factorization categories can be embedded into singularity categories. Let $L$ be a line bundle on $X$, and $W\in \Gamma(X,L)$ be a non-zero-divisor, i.e.~the induced homomorphism $W:\cO_X\to L$ is injective, and denote by $X_0$  the zero scheme of $W$. Denote by $j:X_0\hookrightarrow X$ the closed immersion. Since the direct image $j_*:{\rm D^b}(\coh X_0)\rightarrow{\rm D^b}(\coh X)$ preserves perfect complexes by \cite[Proposition 2.7.(a)]{tt}, it induces an exact functor $$j_{\circ}:{\rm D^{sg}}(X_0)\rightarrow {\rm D^{sg}}(X).$$ As in \cite{orlov4},   the cokernel functor $\Sigma:{\rm MF}(X,L,W)\rightarrow \coh X_0$ defined by $\Sigma(F):={\rm Cok}(\varphi_1^F)$ induces an exact functor
$$\Sigma:{\rm DMF}(X,L,W)\rightarrow{\rm D^{sg}}(X_0).$$

\begin{thm}[{\cite[Theorem 1]{orlov4}}, {\cite[Theorem 2.7]{efi-posi}}]\label{efi-posi}
The functor $\Sigma:{\rm DMF}(X,L,W)\rightarrow{\rm D^{sg}}(X_0)$ is fully faithful, and the essential image of $\Sigma$ is the thick subcategory  consisting of objects $F$ such that $j_{\circ}(F)=0\in{\rm D^{sg}}(X)$. In particular, if $X$ is regular, $\Sigma$ is an equivalence.
\end{thm}

\vspace{2mm}
The following  result is the key motivation for our  definitions of relative singular loci.

\begin{prop}\label{rel sing = supp}
$(1)$ We have an equality of subsets of $X$
$${\rm Sing}(X_0/X)=\bigcup_{F\in {\rm DMF}(X,L,W)}{\rm Supp}(F).$$
$(2)$ We have an equality of subsets of $X$
$${\rm Sing^{loc}}(X_0/X)=\{\,p\in X\,|\,{\rm KMF}(X_p,W_p)\neq0\}$$
\begin{proof}
Since (2) follows from a similar proof of (1), we prove only (1).

If $p\in{\rm Sing}(X_0/X)$, by definition, there exists $A\in{\rm D^{sg}}(X_0)$ such that $A_p\neq0$ in ${\rm D^{sg}}((X_0)_p)$ and  $j_{\circ}(A)=0$ in $\Dsg(X)$. Then by Theorem \ref{efi-posi}, there is an object $F\in{\rm DMF}(X,L,W)$ such that $\Sigma(F)\cong A$. Then we have $A_p\cong\Sigma(F)_p\cong\Sigma_p(F_p)$, where $\Sigma_p:{\rm KMF}(X_p,W_p)\rightarrow {\rm D^{sg}}((X_0)_p)$ is the exact functor defined as above. Since $A_p\neq0$ and $\Sigma_p$ is fully faithful,  $F_p\neq0$ and hence $p\in{\rm Supp}(F)$.

Conversely, if $p\in {\rm Supp}(F)$ for some $F\in{\rm DMF}(X,L,W)$, then $F_p\neq0\in {\rm KMF}(X_p,W_p)$. Since $\Sigma_p:{\rm KMF}(X_p,W_p)\rightarrow {\rm D^{sg}}((X_0)_p)$ is fully faithful, $\Sigma(F)_p\cong\Sigma_p(F_p)\neq0$, and so $\Sigma(F)_p\notin\Perf(X_0)_p$. Furthermore, $j_*(\Sigma(F))\in\Perf X$ by Theorem \ref{efi-posi}. Hence $p\in {\rm Sing}(X_0/X)$.
\end{proof}
\end{prop}

\begin{lem}\label{non zero lemma}

Let $(R,\mathfrak{m})$ be a local ring, and let $W\in R$ be an element.
The category ${\rm KMF}(R,W)$ has a non-zero object if and only if $W\in\mathfrak{m}^2$.
\begin{proof}
Assume that $W\in\m^2$. Then there exist non-units $m_i,n_i\in\m$ for $1\leq i\leq r$ such that $W=\sum_{i=1}^{r}m_in_i$. Set 
$
K_i:=(R\xrightarrow{m_i}R\xrightarrow{n_i}R)\in\KMF(R,m_in_i)
$
and $\cK:=\bigotimes_{i=1}^rK_i\in\KMF(R,W)$. Then we claim that  $\cK\neq0$ in $\KMF(R,W)$. Indeed, if $\cK=0$ in $\KMF(R,W)$, there are morphisms $h_0:\cK_0\to \cK_1$, $h_1:\cK_1\to \cK_0$ such that ${\rm id}_{\cK_0}=\varphi_1^{\cK}h_0+h_1\varphi_0^{\cK}$. Since each $\varphi_i^{\cK}$ is a matrix whose entries are non-units in $R$, the equation implies that $1_{R}\in\m$, which is a contradiction.

For the converse, let $F\in \KMF(R,W)$ is a non-zero object. Since $R$ is local, every locally free modules are free modules. Hence each $\varphi_i^F$ is a $r$-square matrix $(f_{m,n}^i)_{1\leq m,n\leq r}$ in elements in $R$. We claim that $F$ is isomorphic to a matrix factorization $F'$ such that all entries of matrices $\varphi_i^{F'}$ are non-units. Indeed, if there is a unit entry $u\in R$ in the square matrix $\varphi_1^{F}=(f_{m,n}^1)_{1\leq m,n\leq r}$, applying  elementary row/column operations, we may assume that $u=f_{1,1}^1$ and $f^1_{m,1}=f^1_{1,m}=0$ for $m\neq1$. Then we see that $f^0_{m,1}=f^0_{1,m}=0$ for $m\neq1$. Hence the object $\Bigl( R\xrightarrow{u}R\xrightarrow{u^{-1}W}R\Bigr)$ is a direct summand of $F$, but this is isomorphic to the zero object in $\KMF(R,W)$. Since the rank of matrices $\varphi_i^F$ are finite and $F\neq0$, repeating this process, we may assume that  all entries of  $\varphi_i^F$ are non-units, and hence $W\in \m^2$  since $W=\varphi_0^F\varphi_1^F$.
\end{proof}
\end{lem}

The following result is useful to compute the relative singular loci of zero schemes of regular sections of line bundles.

\begin{prop}\label{alg criterion}
Notation is same as above. We have $$\Singloc(X_0/X)=\{\, p\in X\,|\,W_p\in\mathfrak{m}_p^2\,\}$$
\begin{proof}
This follows from Proposition \ref{rel sing = supp}.(2) and Lemma \ref{non zero lemma}.
\end{proof}
\end{prop}

\begin{rem}
 Recall that the \textbf{ critical locus} ${\rm Crit}(\varphi)$ of  a  function $\varphi\in\Gamma(Y,\mathcal{O}_Y)$ on a complex analytic space $(Y,\cO_Y)$ is defined by 
\[
{\rm Crit}(\varphi):=\{ \,p\in Y\,|\, \varphi_p-\varphi(p)\in\mathfrak{m}_p^2\,\},
\]
where $\m_p$ is the maximal ideal of the local ring $\cO_{Y,p}$.
Assume that $X$ is a quasi-projective variety over $\bC$, and let $f\in\Gamma(X,\cO_X)$ be a regular function which is a non-zero-divisor. Denote by $(X^{\rm an},\cO_X^{\rm an})$  the  complex analytic space associated to $X$, and let $f^{\rm an}\in\Gamma(X^{\rm an},\cO_X^{\rm an})$ is the function  associated to $f$. For any $p\in X^{\rm an}$, since the morphism of local rings $\varphi_p:\cO_{X,p}\to\cO_{X,p}^{\rm an}$ is flat \cite{serre}, $\varphi_p$ induces the isomorphism $\m_p/\m_p^2\cong \m_p^{\rm an}/(\m_p^{\rm an})^2$ (see \cite[Theorem 49]{matsumura1}), and so $f_p\in\m_p^2$ if and only if $f^{\rm an}_p\in(\m_p^{\rm an})^2$. Hence
Proposition \ref{alg criterion} implies the following equality of sets;
 \[
\Singloc(f^{-1}(0)/X)\cap X^{\rm an}={\rm Crit}(f^{\rm an})\cap {\rm Zero}(f^{\rm an}),
\]
where ${\rm Zero}(f^{\rm an})$ is the zero locus of $f^{\rm an}$.
\end{rem}

\vspace{3mm}
Recall that  the \textbf{codimension}  of a Noetherian local ring $(R,\m)$ is defined by 
\[
\codim(R):=\embdim(R)-\dim(R),
\]
where $\embdim(R):=\dim_{R/\m}(\m/\m^2)$ is the dimension of $R/\m$-vector space of $\m/\m^2$, which is  called the \textbf{embedding dimension} of $(R,\m)$.  The following result provides a numerical characterization of $\Sing(X_0/X)$.

\begin{prop}
Assume that for any $p\in X_0$ we have $\dim(\cO_{X,p})-\dim(\cO_{{X_0},p})=1$. Then we have the following equality of sets;
\[
\Singloc(X_0/X)=\{\,p\in X_0\,|\, \codim(\cO_{{X_0},p})>\codim(\cO_{X,p})\,\}.
\]
\begin{proof}

Let $p\in X_0$ be a point, and denote by  $\m_p\subset \cO_{X,p}$ and $\n_p\subset \cO_{{X_0},p}$  the maximal ideals of  $\cO_{X,p}$ and $\cO_{{X_0},p}$ respectively. The surjective homomorphism $\pi:\cO_{X,p}\twoheadrightarrow \cO_{{X_0},p}$ induces a surjective map
\[
\overline{\pi}:\m_p/\m_p^2\twoheadrightarrow \n_p/\n_p^2.
\]
Hence we have 
\[
\embdim(\cO_{{X_0},p})\leq \embdim(\cO_{X,p}).
\]
Moreover, since $\dim(\cO_{X,p})-\dim(\cO_{{X_0},p})=1$, we have $\codim(\cO_{{X_0},p})>\codim(\cO_{X,p})$ if and only if  $\embdim(\cO_{{X_0},p})\geq \embdim(\cO_{X,p})$.
Hence, by Proposition \ref{alg criterion}, it is enough to show that 
\[
W_p\in\m_p^2 \hspace{2mm}\Leftrightarrow \hspace{2mm} \overline{\pi} \,\,{\rm is \,\,injective}.
\]

For $(\Rightarrow)$, assume that $W_p\in\m_p^2$, and let $x\in\m_p$ with $\pi(x)\in\n_p^2$. Then, since $\cO_{{X_0},p}\cong\cO_{X,p}/\langle W_p\rangle$, there exists $y\in\cO_{X,p}$ such that 
\[
x+yW_p\in\m_p^2.
\]
Hence $x\in\m_p^2$ by the assumption $W_p\in\m_p^2$. This means that $\overline{\pi}$ is injective.
For $(\Leftarrow)$, assume that $\overline{\pi}$ is injective. Since $p\in X_0$, $W_p\in\m_p$. Since $\pi(W_p)=0$ and $\overline{\pi}$ is injective, we have $W_p\in\m_p^2$.
\end{proof}
\end{prop}

\vspace{2mm}
Next, we show some properties describing  relationships between  relative singular loci and locally relative singular loci.

\begin{prop}\label{singloc and sing}
Let $R$ be a Noetherian ring, and let $W\in R$ is a non-zero-divisor.
Set $X:=\Spec R$ and $X_0:=\Spec (R/W)$.
\begin{itemize}
\item[$(1)$] We have 
\[
\Sing (X_0/X)\cap\Max (R/W)=\Sing^{\rm loc}(X_0/X)\cap\Max (R/W).
\]
\item[$(2)$] Let $p\in X_0$ be a point with $p\notin\Singloc(X_0/X)$. Let $q\in X_0$ be a point  and denote by $\overline{\{q\}}$ the closure of $q$. If $p\in \overline{\{q\}}$, then $q\notin\Sing(X_0/X)$.

\end{itemize}
\begin{proof}
(1) The inclusion $\Sing (X_0/X)\cap\Max (R/W)\subseteq\Sing^{\rm loc}(X_0/X)\cap\Max (R/W)$ follows from Proposition \ref{rel sing = sing}.  To show the opposite inclusion, let $\m\in\Sing^{\rm loc}(X_0/X)\cap\Max (R/W)$ be a maximal ideal of $R/W$ which is contained in $\Singloc(X_0/X)$, and let $\overline{\m}\in \Max R$ be the maximal ideal of $R$ such that $i(\m)=\overline{\m}$, where $i:X_0\hookrightarrow X$ is the closed immersion. By Proposition \ref{alg criterion}, there exist elements $m_i,n_i,r\in R$ ($1\leq i\leq r$) such that  $m_i\in\overline{\m}$, $n_i\in\overline{\m}$, $r\notin\overline{\m}$, and $rW=\sum_{i=1}^{r}m_in_i$. Since $r\notin\overline{\m}$ and $\overline{\m}$ is maximal, we have $\langle r\rangle +\overline{\m}=R$, and so there exists an element $a\in R$ such that $1-ar\in\overline{\m}$.  Then we have 
\[W=(1-ar)W+\sum_{i=1}^{r}am_in_i.\]
Consider the following matrix factorizations  
\begin{eqnarray*}
K_0:=\Bigl(R\xrightarrow{W}R\xrightarrow{1-ar}R\Bigr)\hspace{2mm}\in&\KMF(R,(1-ar)W)\\
K_i:=\Bigl(R\xrightarrow{am_i}R\xrightarrow{n_i}R\Bigr)\hspace{2mm}\in&\hspace{-6mm}\KMF(R,am_in_i)\\
\cK:=\bigotimes_{i=0}^{r}K_i\hspace{2mm}\in&\hspace{-11mm}\KMF(R,W).
\end{eqnarray*}
Then, by Lemma \ref{koszul supp} and Lemma \ref{supp tens}, we see that  $\overline{\m}\in\bigcap_{i=0}^{r}\Supp (K_i)=\Supp (\cK)$. Proposition \ref{rel sing = supp}.(1) implies that  $\m\in\Sing(X_0/X)$.

(2) Assume that $q\in\Sing(X_0/X)$. Since the relative singular locus $\Sing(X_0/X)$ is a union of closed subsets by Proposition \ref{rel sing = supp}.(1), we have $\overline{\{q\}}\subseteq\Sing(X_0/X)$.  Since $p\in\overline{\{q\}}$ and $\Sing(X_0/X)\subseteq\Singloc(X_0/X)$, $p\in\Singloc(X_0/X)$, which contradicts to the assumption of $p$.
\end{proof}
\end{prop}

\vspace{2mm}
At the end of this section, we compute examples of relative singular loci using the above results.

\begin{exa}\label{example}
We give two examples of the relative singular loci which are not equal to the usual singular loci.

$(1)$ Let $R:=\bC[x,y]/\langle x^n\rangle$ for $n>1$, and let $W:=\overline{y}\in R$. Set $X:=\Spec R$ and $X_0:=\Spec (R/W)$. Although $\Sing(X_0)=\{{\rm pt}\}\neq\emptyset$, by Proposition \ref{alg criterion}, we have
\[\Sing(X_0/X)=\Sing^{\rm loc}(X_0 /X)=\emptyset.\]

$(2)$ Let $R:=\bC[x,y,z,w]/\langle xy-zw\rangle$, and let $W:=\overline{w}\in R$. Set $X:=\Spec R$ and $X_0:=\Spec (R/W)$.  Then we have \[\Sing(X_0)=\{\,\langle \overline{x},\overline{y}\rangle, \langle \overline{x},\overline{y},\overline{z}-a\rangle\,|\, a\in\bC\,\}.\]
By Proposition \ref{alg criterion}, we see that 
\[\Singloc(X_0 /X)=\Sing(X_0)\setminus\{{\langle \overline{x},\overline{y},\overline{z}\rangle}\},\]
and, by Proposition \ref{singloc and sing}, we have 
\[
\Sing(X_0/X)=\Singloc(X_0/X)\setminus \{{\langle \overline{x},\overline{y}\rangle}\}.
\]
In this example, all kinds of singular loci are different;
\[\Sing(X_0/X)\subsetneq\Singloc(X_0/X)\subsetneq\Sing(X_0).\]
\end{exa}

\vspace{2mm}
\section{Tensor nilpotence properties}

In this section, we prove the tensor nilpotent properties, which will be necessary for our main result. 
The properties are analogous to \cite[Theorem 3.6, 3.8]{thomason}, and the strategy of the proof  is similar to loc.~cit.
Let $X$ be a  Noetherian scheme, and let $W\in\Gamma(X,L)$ be a global section of a line bundle $L$ on $X$.

\subsection{Mayer-Vietoris sequence}

We provide a Mayer-Vietoris  sequence for factorizations for the proof of the tensor nilpotence properties in the next section.
For an open immersion $i:U\hookrightarrow X$,  consider an induced LG model $(U,L|_U,W|_U)$. Let $i_!:\Mod\mathcal{O}_U\rightarrow \Mod\mathcal{O}_X$ be the extension by zero. For an object $F=\Bigl(F_1\xrightarrow{\varphi_1^F} F_0\xrightarrow{\varphi_0^F} F_1\otimes L|_U\Bigr)\in{\rm Sh}(U,L|_U,W|_U)$, we define an object $i_!(F)\in{\rm Sh}(X,L,W)$ by
$$i_!(F):=\Bigl(i_!(F_1)\xrightarrow{i_!(\varphi_1^F)} i_!(F_0)\xrightarrow{\sigma\circ i_!(\varphi_0^F)} i_!(F_1)\otimes L\Bigr),$$
where $\sigma:i_!(F_1\otimes L|_U)\xrightarrow{\sim} i_!(F_1)\otimes L$ is a natural isomorphism. This defines an exact functor $$i_!:{\rm Sh}(U,L|_U,W|_U)\rightarrow{\rm Sh}(X,L,W).$$
Similarly, the inverse image functor $i^*:\Mod\mathcal{O}_X\rightarrow \Mod\mathcal{O}_U$ defines an exact functor 
$$i^*:{\rm Sh}(X,L,W)\rightarrow{\rm Sh}(U,L|_U,W|_U).$$ 
These functors induce an exact functors between homotopy categories; 
$$i_!:{\rm KSh}(U,L|_U,W|_U)\rightarrow{\rm KSh}(X,L,W)$$ 
$$i^*:{\rm KSh}(X,L,W)\rightarrow{\rm KSh}(U,L|_U,W|_U).$$
Since both of $i_!$ and $i^*$ are exact functors and preserve arbitrary direct sums, these  functors defines exact functors between coderived categories
$$i_!:{\rm D^{co}Sh}(U,L|_U,W|_U)\rightarrow{\rm D^{co}Sh}(X,L,W)$$
$$i^*:{\rm D^{co}Sh}(X,L,W)\rightarrow{\rm D^{co}Sh}(U,L|_U,W|_U).$$
Adjunction $i_!\dashv i^*$ of functors between $\Mod\mathcal{O}_U$ and $\Mod\mathcal{O}_X$ naturally induces adjunction $i_!\dashv i^*$ of functors between coderived categories ${\rm D^{co}Sh}(U,L|_U,W|_U)$ and ${\rm D^{co}Sh}(X,L,W)$.

\begin{lem}[Mayer-Vietoris]\label{mayer}
Let $U_1$ and $U_2$ be open subschemas of  $X$, and suppose $X=U_1\cup U_2$.
Denote by $U_{1,2}$ the intersection $U_1\cap U_2$. Let $i_l:U_l\rightarrow X$ and $i_{1,2}:U_{1,2}\rightarrow X$ be open immersions. Then, for any objects $E, F\in {\rm D^{co}Sh}(X,L,W)$ and for any $k\in\mathbb{Z}$,  we have the following long exact sequence:


\begin{tikzcd}
  {\rm Hom}_{U_{1,2}}(i_{1,2}^*E,i_{1,2}^*F[k-1])  \rar{\delta_{k-1}} & {\rm Hom}_X(E,F[k])  \rar
             \ar[draw=none]{d}[name=X, anchor=center]{} 
             & \bigoplus_{l=1}^{2}{\rm Hom}_{U_l}(i_{l}^*E,i_{l}^*F[k]) \ar[rounded corners,
            to path={ -- ([xshift=2ex]\tikztostart.east)
                      |- (X.center) \tikztonodes
                      -| ([xshift=-2ex]\tikztotarget.west)
                      -- (\tikztotarget)}]{dll}[at end]{} \\
  {\rm Hom}_{U_{1,2}}(i_{1,2}^*E,i_{1,2}^*F[k]) 
  \rar{\delta_k} &   {\rm Hom}_X(E,F[k+1]) \rar &   \bigoplus_{l=1}^{2}{\rm Hom}_{U_l}(i_{l}^*E,i_{l}^*F[k+1])
   \end{tikzcd}\vspace{1mm}\\
    where ${\rm Hom}_{(-)}$ denotes the set of morphisms in ${\rm D^{co}Sh}(-,L|_{(-)},W|_{(-)})$.

\begin{proof}
  For each components $E_j$ ($j=0,1$), we have the following exact sequence in $\Mod\mathcal{O}_X$:
$$0\rightarrow {i_{1,2}}_!i_{1,2}^*E_j\rightarrow{i_{1}}_!i_{1}^*E_j\oplus{i_{2}}_!i_{2}^*E_j\rightarrow E_j\rightarrow 0.$$
These exact sequences induce the following exact sequence in the abelian category ${\rm Sh}(X,L,W)$:
$$0\rightarrow {i_{1,2}}_!i_{1,2}^*E\rightarrow{i_{1}}_!i_{1}^*E\oplus{i_{2}}_!i_{2}^*E\rightarrow E\rightarrow 0.$$
By the same argument as in \cite[Lemma 2.7.(a)]{ls}, the above exact sequence gives rise to a triangle 
$$(\ast)\hspace{7mm}{i_{1,2}}_!i_{1,2}^*E\rightarrow{i_{1}}_!i_{1}^*E\oplus{i_{2}}_!i_{2}^*E\rightarrow E\rightarrow{i_{1,2}}_!i_{1,2}^*E[1]$$
in ${\rm D^{co}Sh}(X,L,W)$. By applying ${\rm Hom}_{{\rm D^{co}Sh}(X,L,W)}(-,F[k])$ to the triangle and using the adjunction ${i_{(-)}}_!\dashv {i_{(-)}}^*$ for open immersions $i_{(-)}$, we obtain the result.
\end{proof}  
\end{lem}

\subsection{Tensor nilpotence properties}

The following lemma is an analogy of \cite[Theorem 3.6]{thomason}, and we show it by  a similar argument in the proof of loc. cit.

\begin{lem}\label{tens nilp1}
Let $f:E\rightarrow F$ be a morphism in ${\rm DMF}(X,L,W)$. If $f\otimes k(p)=0$ in ${\rm DMF}(k(p),W\otimes k(p))$ for any $p\in X$, then there is an integer $n$ such that $f^{\otimes n}=0$ in ${\rm DMF}(X,L,nW)$.\begin{proof}
The proof will be divided into 5 steps.

\vspace{1mm}
{\sf Step~1.} In the first step, we will reduce to the case that $X$ is affine. Since $X$ is Noetherian, we have a finite number of open affine covering $\bigcup_{i=1}^k{U_i}$ of $X$, where $U_i=\Spec R_i$. Set $L_i:=L|_{U_i}$ and $W_i:=W|_{U_i}$. We will show that, if there are positive integers $n_i$ such that $(f|_{U_i})^{\otimes n_i}=0\in{\rm KMF}(R_i,L_i,n_iW_i)$ for every $i$, there is a positive integer $n$ such that $f^{\otimes n}=0\in {\rm DMF}(X,L,nW)$. We will do this by induction on $k$. If $k=1$, then $f^{\otimes n_1}=0$. For $k>1$, suppose that the result is true for $ k-1$. Set $V:=\bigcup_{i=2}^kU_i$. Then by induction hypothesis, there exists $n'$ such that $(f|_V)^{\otimes n'}=0$ in ${\rm DMF}(V,L|_V,n'W|_V)$. If we set $m:={\rm max}\{n_1,n'\}$, we have $f^{\otimes m}|_V=0$ and $f^{\otimes m}|_{U_1}=0$. 
By Lemma \ref{mayer}, there exists a morphism $g:j^*E^{\otimes m}[1]\rightarrow j^*F^{\otimes m}$ in ${\rm D^{co}Sh}(U_1\cap V,L|_{U_1\cap V},W|_{U_1\cap V})$ such that the following diagram in ${\rm D^{co}Sh}(X,L,W)$ is commutative:
\[\xymatrix{
E^{\otimes m}\ar[r]^{f^{\otimes m}}\ar[d]_{\delta}&F^{\otimes m}\\
j_!(j^*E^{\otimes m}[1])\ar[r]^{j_!(g)}&j_!j^*F^{\otimes m}\ar[u]_{\varepsilon}
}\]
where $j:U_1\cap V\rightarrow X$ is the open immersion, $\delta$ is the morphism corresponding to the third morphism in the above triangle $(\ast)$ in the proof of Lemma \ref{mayer}, and $\varepsilon$ is the adjunction morphism of $j_!\dashv j^*$. Since $f^{\otimes m}\otimes j_!(g)$ is identified with $j_!(j^*(f^{\otimes m})\otimes g)$ via natural isomorphisms, we have $f^{\otimes m}\otimes j_!(g)=0$. Hence for $n:=2m$, it follows that $f^{\otimes n}=0$  in ${\rm D^{co}Sh}(X,L,nW)$, as $f^{\otimes n}$ factors through $f^{\otimes m}\otimes j_!(g)$. By Lemma \ref{ff}, we obtain $f^{\otimes n}=0$ in ${\rm DMF}(X,L,nW)$.

\vspace{2mm}

{\sf Step 2.} By the above step, we may assume that  $X=\Spec R$ is affine and $L$ is an invertible $R$-module. Recall that ${\rm DMF}(R,L,W)\cong {\rm KMF}(R,L,W)$ by Proposition \ref{affine}. In this step, we reduce to the case when $W=0$ and $E=(0\rightarrow R\rightarrow 0)$. We have the following adjunction
$$\Phi:{\rm Hom}_{{\rm KMF}(R,L,W)}(E,F)\xrightarrow{\sim}{\rm Hom}_{{\rm KMF}(R,L,0)}(R,E^{\vee}\otimes F).$$
Via natural isomorphisms, we can identify $\Phi(f^{\otimes n})$ with $\Phi(f)^{\otimes n}$ and $\Phi(f\otimes k(p))$ with $\Phi(f^{\otimes n})\otimes k(p)$ respectively. Therefore, we may assume that $W=0$ and $E=R$.
\vspace{2mm}

{\sf Step 3.}  Since the components $F_i$ of $F$ and $L$ are locally free, by the above first step, shrinking $X=\Spec R$ if necessary, we may assume that $X=\Spec R$ is affine scheme such that  each $F_i$ is free $R$-module and $L\cong R$. Furthermore, by the second step, we may assume that $E=R$ and $W=0$.
Hence, since the natural isomorphism $\Hom_{\MF(R,0)}(R,F)\cong {\rm Ker}(\varphi_0^F)$ induces the isomorphism ${\rm Hom}_{{\rm KMF}(R,0)}(R,F)\cong {\sf H}_0(F)$, we only need to show that if  $f\in {\sf H}_0(F)$ satisfies $f\otimes k(p)=0\in{\sf H}_0(F\otimes k(p))$ for any $p\in \Spec R$, then $f^{\otimes n}=0\in{\sf H}_0(F^{\otimes n})$ for some $n>0$.

In this step, we reduce to the case when the Noetherian ring $R$ is of finite Krull dimension. Since the components $F_i$ of $F$ are free $R$-modules, the morphisms $\varphi_i^F$ can be represented by a matrices whose entries are elements in $R$. Let $\{R_{\alpha}\}_{\alpha\in\cA}$ be the family of all subrings $R_{\alpha}\subset R$ of  $R$ such that $\dim R_{\alpha}<\infty$ and $R_{\alpha}$ contains all entries of matrices $\varphi_1$ and $\varphi_0$. 
Then for any $\alpha\in\cA$, there is the natural object $F_\alpha\in\KMF(R_{\alpha},0)$ such that its components $(F_\alpha)_i$  are free $R_{\alpha}$-modules and $\pi_{\alpha}^*(F_{\alpha})\cong F$ in the additive category $\MF(R,0)$, where $\pi_{\alpha}:\Spec R\rightarrow \Spec R_{\alpha}$ is the morphism induced by the inclusion $R_{\alpha}\subset R$.
Let \[F^{\text{\tiny{\textbullet}}}:=(\cdots0\rightarrow F_1\xrightarrow{\varphi_1^F} F_0\xrightarrow{\varphi_0^F}F_1\rightarrow0\cdots )\]
\[F^{\text{\tiny{\textbullet}}}_{\alpha}:=(\cdots0\rightarrow (F_{\alpha})_1\xrightarrow{\varphi_1^{F_{\alpha}}} (F_{\alpha})_0\xrightarrow{\varphi_0^{F_{\alpha}}}(F_{\alpha})_1\rightarrow0\cdots )\]
be the complexes of free modules such that  the term $F_0$ and $(F_{\alpha})_0$ are of degree $0$. Then we have $H^0(F^{\text{\tiny{\textbullet}}})={\sf H}_0(F)$ and $H^0(F^{\text{\tiny{\textbullet}}}\otimes k(p))={\sf H}_0(F\otimes k(p))$ for any $p\in \Spec R$.
Since  $R$ is the direct colimit of the system $\{R_{\alpha}\}_{\alpha\in\cA}$; $R=\colim R_{\alpha}$, 
by the same argument as in the  step (3.6.4) in the proof of \cite[Theirem 3.6]{thomason}, if $f\otimes k(p)=0$ in ${\sf H}_0(F\otimes k(p))$ for any $p\in X$, there exist $\beta\in\cA$ and an element $f_{\beta}\in H^0(F^{\text{\tiny{\textbullet}}}_{\beta})={\sf H}_0(F_{\beta})$ such that $\pi^*_{\beta}(f_{\beta})=f$ and
$f_{\beta}\otimes k(p)=0$ in $H^0(F^{\text{\tiny{\textbullet}}}_{\beta}\otimes k(p))={\sf H}_0(F_{\beta}\otimes k(p))$ for any $p\in\Spec R_{\beta}$. Therefore, if the assertion is true for $X=\Spec R_{\beta}$, there exists $n>0$ such that $f_{\beta}^{\otimes n}=0$ in ${\sf H}_0(F_{\beta}^{\otimes n})$, and  then $f^{\otimes n}=\pi^*_{\beta}(f_{\beta}^{\otimes n})=0$ in ${\sf H}_0(F^{\otimes n})$.
This completes the reduction to the case when the ring $R$ is  of finite Krull dimension.

\vspace{2mm}

{\sf Step 4.}
In this step, we reduce to the case when the Noetherian ring $R$ of finite Krull dimension is reduced.
Let $\mathfrak{N}\subset R$ be the ideal of nilpotent elements in $R$, and denote by $h:\Spec R/\mathfrak{N}\rightarrow \Spec R$ be the closed immersion. If  $f\otimes k(p)=0$ for any $p\in \Spec R$, then $h^*(f)\otimes k(q)=0$ for any $q\in\Spec R/\mathfrak{N}$. Hence, for the reduction, we  claim that, if $h^*(f)^{\otimes n'}=0\in{\sf H}_0((h^*F)^{\otimes n'})$ for some $n'>0$, then there exists $n>0$ such that $f^{\otimes n}=0\in {\sf H}_0(F^{\otimes n})$.  Assume that $h^*(f)^{\otimes n'}=0\in{\sf H}_0((h^*F)^{\otimes n'})$. Since $(h^*F)^{\otimes n'}\cong h^*(F^{\otimes n'})=F^{\otimes n'}\otimes _{R}R/\mathfrak{N}$, the assumption implies that there exist  elements $x\in (F^{\otimes n'})_1$ and $y\in\mathfrak{N}\bigl((F^{\otimes n'})_0\bigr)\subset (F^{\otimes n'})_0$ such that $f^{\otimes n'}=\varphi_{1}^{F^{\otimes n'}}(x)+y$ in $(F^{\otimes n'})_0$. Since $y\in\mathfrak{N}\bigl((F^{\otimes n'})_0\bigr)$, there is a positive integer $m$ such that $y^{\otimes m}=0$ in the free $R$-module $\bigl((F^{\otimes n'})_0\bigr)^{\otimes m}$. Therefore, it is enough to show the following claim:

 \vspace{2mm}
 Let $S$ be a ring, and let  $E\in {\rm MF}(S,0)$ be an object such that  its components $E_i$ are free $S$-modules. For $e\in {\rm Ker}(\varphi_0^E)\subset E_0$, suppose that  there are elements $u\in E_1$ and $v\in E_0$ such that $e=\varphi_1^{E}(u)+v$ and $v^{\otimes n}=0$ in the $S$-free module $(E_0)^{\otimes n}$ for some $n>0$. Then, considering $e^{\otimes n}$ as an element in $(E^{\otimes n})_0$ via the natural split mono $(E_0)^{\otimes n}\hookrightarrow (E^{\otimes n})_0$, there is an element $w\in (E^{\otimes n})_1$ such that $\varphi_1^{E^{\otimes n}}(w)=e^{\otimes n}$ in $(E^{\otimes n})_0$. In particular, $e^{\otimes n}=0\in {\sf H}_0(E^{\otimes n})$.

\vspace{2mm}
 The element $e^{\otimes n}=(\varphi_1^{E}(u)+v)^{\otimes n}\in (E_0)^{\otimes n}$ can be decomposed into the following form
$$e^{\otimes n}=\varphi_1^E(u)\otimes w_{n-1}+v\otimes\varphi_1^E(u)\otimes w_{n-2}+\cdots +v^{\otimes n-2}\otimes \varphi_1^E(u)\otimes w_1+ v^{\otimes n-1}\otimes \varphi_1^E(u),$$
where $w_i$ is an element in $(E_0)^{\otimes i}$. For an ordered sequence  $(i_1,i_2,...,i_n)$ of $i_k\in\{0,1\}$, set $E_{(i_1,i_2,...,i_n)}:=E_{i_1}\otimes E_{i_2}\otimes \cdots\otimes E_{i_n}$, and set $$\widetilde{E}:=E_{(1,0,0,...,0)}\oplus E_{(0,1,0,...,0)}\oplus\cdots\oplus E_{(0,...,0,1)}.$$
Let $\widetilde{w}:=(u\otimes w_{n-1})\oplus (v\otimes u\otimes w_{n-2})\oplus\cdots \oplus (v^{\otimes n-2}\otimes u\otimes w_1)\oplus (v^{\otimes n-1}\otimes u)\in \widetilde{E}$, and let $w:=\iota(\widetilde{w})\in (E^{\otimes n})_1$ be the image of $\widetilde{w}$ under the natural split mono $\iota:\widetilde{E}\hookrightarrow (E^{\otimes n})_1$. Since $\varphi_0^E(\varphi_1^E(u))=0$, $\varphi_0^E(v)=\varphi_0^E(e-\varphi_1^E(u))=0$,  and each $w_i\in(E_0)^{\otimes i}$ is a summation of elements of the form $a_1\otimes a_2\otimes \cdots\otimes a_i$ where $a_k$ are either $\varphi_1^E(u)$ or $v$, we obtain an equality $\varphi_1^{E^{\otimes n}}(w)=e^{\otimes n}$ in $(E^{\otimes n})_0$. This completes the proof of the claim.

\vspace{2mm}
{\sf Step 5.}
Now we may assume that $X=\Spec R$ is reduced affine scheme of finite Krull dimension, the components $F_i$ of $F$ are free $R$-modules, $L\cong R$, $W=0$, and $E=R\in {\rm KMF}(R,0)$. In this step, we finish the proof by induction on the Krull dimension $d:=\dim R$ of $R$. 

If $d=0$, then $R\cong\bigoplus_{p\in X} k(p)$ and  ${\sf H}_0(F)\cong \bigoplus_{p\in X}{\sf H}_0(F\otimes k(p))$. Hence, if $f\otimes k(p)=0\in{\sf H}_0(F\otimes k(p))$ for any $p\in X$, then $f=0$ in ${\sf H}_0(F)$. 

Consider a case when $d>0$, and assume that the result holds for Noetherian rings of dimension less than $d$. Denote by $\Min R$  the finite set of all prime ideals of $R$ of height zero. Then the product $\prod_{p\in\Min R}k(p)$ of residue fields is isomorphic to the localization $S^{-1}R$ for the set $S$ of all non zero-divisors in $R$, as the residue fields $k(p)$ is equal to the local ring $R_{p}$ for any $p\in\Min R$. By hypothesis,  $f\otimes k(p)=0\in{\sf H}_0(F\otimes k(p))$ for any $p\in\Min R$, hence $f\otimes_R S^{-1}R=0\in {\sf H}_0(F\otimes _R S^{-1}R)$. This means that, for a representative $f\in F_0$ of the equivalence class $f\in{\sf H}_0(F)$, there exist elements $y\in F_1$ and $s\in S$ such that $sf=\varphi_1^F(y)$ in $F_0$.
Set $K_s:=(R\xrightarrow{s}R\xrightarrow{0}R)\in{\rm MF}(R,0)$ and let $\gamma:=(y,f):K_s\rightarrow F$ be the  morphism defined as the pair of morphisms $y:R\rightarrow F_1$ and $f:R\rightarrow F_0$ of $R$-modules. Denote by $i:\Spec R/s\rightarrow \Spec R$ the natural closed immersion. The canonical quotient $R\rightarrow i_*(R/s)$ naturally defines morphisms $\delta:K_s\rightarrow (0\rightarrow i_*(R/s)\rightarrow 0)$ and $\alpha:(0\rightarrow R\rightarrow 0)\rightarrow (0\rightarrow i_*(R/s)\rightarrow 0)$ in ${\rm coh}(R,0)$. Since $s:R\rightarrow R$ is injective, we have an exact sequence 
$$0\rightarrow \Bigl(R\xrightarrow{\rm id}R\xrightarrow{0}R\Bigr)\xrightarrow{({\rm id},s)} K_s\xrightarrow{\delta}\Bigr(0\rightarrow i_*(R/s)\rightarrow 0\Bigr)\rightarrow 0 $$
in ${\rm coh}(R,0)$. Since $(R\xrightarrow{\rm id}R\xrightarrow{0}R)$ is zero in ${\rm Dcoh}(R,0)$, $\delta$ is an isomorphism in ${\rm Dcoh}(R,0)$  by \cite[Lemma 2.7.(a)]{ls}. Set $\beta:(0\rightarrow i_*(R/s)\rightarrow 0)\rightarrow F$ be the composition $\gamma\circ\delta^{-1}$ in ${\rm Dcoh}(R,0)$. Then the composition $\beta\circ\alpha:(0\rightarrow R\rightarrow 0)\rightarrow F$ is equal to $f:(0\rightarrow R\rightarrow 0)\rightarrow F$ in ${\rm Dcoh}(R,0)$, since $\alpha =\delta\circ\iota$ and $f=\delta \circ \iota$ in ${\rm coh}(R,0)$, where $\iota: (0\rightarrow R\rightarrow 0)\rightarrow K_s$ is the morphism such that $\iota_1=0$ and $\iota_0={\rm id}_R$. Hence,  for any $n>0$, we have the following commutative diagram in ${\rm Dcoh}(R,0)$:
\[\xymatrix{
R\ar[r]^-{\sim}\ar[dd]_-{f^{\otimes n+1}}&R^{\otimes n}\otimes R\ar[dd]_{f^{\otimes n}\otimes f}\ar[dr]^{f^{\otimes n}\otimes 1}\ar[rr]^{1\otimes \alpha}&&R^{\otimes n}\otimes i_*(R/s)\ar[dd]^-{f^{\otimes n}\otimes 1}\ar[r]^-{\sim}&i_*(R/s^{\otimes n})\ar[dd]^{i_*(i^*f^{\otimes n})}\\
&&F^{\otimes n}\otimes R\ar[dl]_{1\otimes f}\ar[dr]^{1\otimes \alpha}&&\\
F^{\otimes n+1}\ar[r]^-{\sim}&F^{\otimes n}\otimes F&&F^{\otimes n}\otimes i_*(R/s)\ar[ll]^{1\otimes \beta}\ar[r]^-{\sim}&i_*(i^*F^{\otimes n})
}\]
where $R=(0\rightarrow R\rightarrow0)$ and $i_*(R/s)=(0\rightarrow i_*(R/s)\rightarrow0)$ by abuse of notation. Since $\dim R/s<\dim R$ and $i^*f\otimes_{R/s} k(p)=f\otimes_R k(p)=0$ in ${\sf H}_0(i^*F\otimes k(p))$, by the induction hypothesis, there exists $m>0$ such that $i^*f^{\otimes m}=0$ in ${\rm KMF}(R,0)$, in particular $i^*f^{\otimes m}=0$ in ${\rm Dcoh}(R,0)$. By the above commutative diagram, we see that  $f^{\otimes m+1}$ factors through $i_*(i^*f^{\otimes m})$ in ${\rm Dcoh}(R,0)$, and hence $f^{\otimes m+1}=0$ in ${\rm Dcoh}(R,0)$. Then, by Lemma \ref{ff}.(5), $f^{\otimes m+1}=0$ in ${\rm KMF}(R,0)$, and  this completes the proof.
\end{proof}
\end{lem}

The following lemma is a consequence of the above lemma.

\begin{lem}\label{tens nilp2}
Let $a:E\rightarrow F$ be a morphism in ${\rm DMF}(X,L,0)$, and let $G\in{\rm DMF}(X,L,W)$ be an object.
If  $a\otimes k(p)=0$ in ${\rm DMF}(k(p),0)$ for all $p\in {\rm Supp}(G)$, then there is an integer $n>0$ such that $G\otimes (a^{\otimes n})=0$ in ${\rm DMF}(X,L,W)$. 
\begin{proof}
Since $(G\otimes a)\otimes k(p)=0$ in ${\rm DMF}(k(p),W\otimes k(p))$  for any $p\in X$, by Lemma \ref{tens nilp1}, there is a positive integer $n>0$ such that $(G\otimes a)^{\otimes n}=G^{\otimes n}\otimes a^{\otimes n}=0$ in ${\rm DMF}(X,L,nW)$. Hence $((G^{\vee})^{\otimes n-1}\otimes G^{\otimes n})\otimes a^{\otimes n}=0$, and so it suffices to show that $G\otimes a^{\otimes n}$ is a retract of $((G^{\vee})^{\otimes n-1}\otimes G^{\otimes n})\otimes a^{\otimes n}$ in ${\rm DMF}(X,L,W)$.

We will show it by proving that  $G$ is a direct summand of $(G^{\vee})^{\otimes n-1}\otimes G^{\otimes n}$ in ${\rm DMF}(X,L,W)$ by induction on $n$. The $n=1$ case is trivial. For $n\geq 2$, assume that $G$ is a direct summand of  $(G^{\vee})^{\otimes n-2}\otimes G^{\otimes n-1}$. It is enough to show that $(G^{\vee})^{\otimes n-2}\otimes G^{\otimes n-1}$ is a direct summand of $(G^{\vee})^{\otimes n-1}\otimes G^{\otimes n}$. But for $n\geq3$ case, this follows from $n=2$ case by tensoring $(G^{\vee})^{\otimes n-3}\otimes G^{\otimes n-2}$. Therefore it suffices to prove that $G$ is a direct summand of $G^{\vee}\otimes G^{\otimes 2}$. The tensor product $(-)\otimes G:{\rm DMF}(X,L,0)\rightarrow {\rm DMF}(X,L,W)$ is left adjoint to $(-)\otimes G^{\vee}:{\rm DMF}(X,L,W)\rightarrow{\rm DMF}(X,L,0)$. Let $\eta:{\rm id}\rightarrow(-)\otimes G\otimes G^{\vee}$ be its adjunction morphism. Then the functor morphism 
$$(-)\otimes G\xrightarrow{\eta\otimes G} (-)\otimes G\otimes G^{\vee}\otimes G$$
is a split mono. Evaluating $\mathcal{O}_X\in {\rm DMF}(X,L,0)$, we see that $G$ is a direct summand of $G^{\vee}\otimes G^{\otimes 2}$.

\end{proof}
\end{lem}

\vspace{2mm}
\section{Classification of thick subcategories of ${\rm DMF}(X,L,W)$}

In this section, we prove our main result. Let $X$ be a separated Noetherian scheme, and let $W\in\Gamma(X,L)$ be any section of a line bundle $L$ on $X$. At first,  following \cite{tt}, we recall the definition   of ample families of line bundles.

\begin{dfn}
We say a quasi-compact and quasi-separated scheme $S$ \textbf{has an ample family of line bundles} if there exists a family $\cL:=\{L_{\alpha}\}_{\alpha\in\cA}$  of line bundles on $S$ such that the family $\{S_f\,|\,f\in \Gamma(S,L_{\alpha}^{\otimes n}), L_{\alpha}\in\cL, n>0\}_{}$ of open subsets form an open basis of $S$, where $S_f:=\{\,p\in S\,|\, f(p)\neq0\}$.
\end{dfn}

\begin{rem}

(1) Any scheme with an ample line bundle has an ample family of line bundles. In particular, any affine scheme has an ample family of line bundles. Any separated regular Noetherian scheme has an ample family of line bundles. See \cite[Example 2.1.2]{tt} for more examples of schemes with ample families of line bundles.\\
(2) If $S$ has an ample family of line bundles, then $S$ satisfies the resolution property by  \cite[Lemma 2.1.3]{tt}.
\end{rem}

\begin{prop}\label{support exist}
 Let $Z\subseteq X$ be a closed subset of $X$. Assume that $X$ has an ample family of line bundles. Then, the following holds.
\begin{itemize}
 \item[$(1)$] There exists a matrix factorization $\cK\in {\rm DMF}(X,L,0)$ such that ${\rm Supp}(\cK)=Z$.
 
 \item[$(2)$] If $W$ is a non-zero-divisor and $Z$ is contained in ${\rm Sing}(X_0/X)$, then there exists a matrix factorization $F\in{\rm DMF}(X,L,W)$ such that ${\rm Supp}(F)=Z$.
 \end{itemize}

 \begin{proof}
 $(1)$ Since X has an ample family of line bundles, there are finitely many sections $f_i\in\Gamma(X,L_i)$ of line bundles $L_i$ ($1\leq i\leq l$) such that $Z=\bigcap_{i=1}^{l} Z(f_i)$ as closed subsets of $X$, where $Z(f_i)$ is the zero scheme of $f_i$. Let $K_i$ be the object in ${\rm DMF}(X,L,0)$ of the following form
 $$K_i:=\Bigl(\mathcal{O}_X\xrightarrow{f_i}L_i\xrightarrow{0}L\Bigr).$$
 By Proposition \ref{supp local}.(1) and Lemma \ref{koszul supp}, we have ${\rm Supp}(K_i)=Z(f_i)$.  If we set $\cK:=\bigotimes_{i=1}^{l} K_i$,  then $\cK\in{\rm DMF}(X,L,0)$ and ${\rm Supp}(\cK)=\bigcap{\rm Supp}(K_i)=Z$.

 $(2)$ Since $X$ is Noetherian, we can decompose $Z$ into finitely many irreducible components $Z=\bigcup_{I=1}^{r} Z_i$. Then, for each $1\leq i\leq r$, there is a unique generic point $p_i\in Z_i$ of $Z_i$. By Proposition \ref{rel sing = supp}, there exists a matrix factorization $E_i\in{\rm DMF}(X,L,W)$ such that $p_i\in{\rm Supp}(E_i)$. Then $Z_i\subseteq {\rm Supp}(E_i)$, since $p_i$ is a generic point of $Z_i$. By $(1)$, there are matrix factorizations $\cK_i\in{\rm DMF}(X,L,0)$, for $1\leq i\leq r$, such that ${\rm Supp}(\cK_i)=Z_i$.
 Note that Lemma \ref{supp tens} implies the equality ${\rm Supp}(\cK_i\otimes E_i)={\rm Supp}(\cK_i)\cap{\rm Supp}(E_i)$.
If we set  $F:=\bigoplus_{i=1}^{r}({\cK_i\otimes E_i})$, then $F\in{\rm DMF}(X,L,W)$ and we have 
$${\rm Supp}(F)=\bigcup_{i=1}^{r}{\rm Supp}(\cK_i\otimes E_i)=\bigcup_{i=1}^{r}({\rm Supp}(\cK_i)\cap{\rm Supp}(E_i))=Z.$$
This completes the proof.
  \end{proof}
\end{prop}

\begin{dfn}
Let $\cT\subset \DMF(X,L,W)$ be a triangulated  full subcategory. We say that $\cT$ is \textbf{$\otimes$-submodule} if it is closed under tensor action of $\DMF(X,L,0)$, i.e.~for any $F\in\DMF(X,L,0)$ and any $T\in\cT$, we have $F\otimes T\in \cT$.   For an object $F\in \DMF(X,L,W)$, we denote by 
\[
\langle F\rangle^{\otimes}\subset \DMF(X,L,W)
\]
the smallest thick $\otimes$-submodule containing $F$.
\end{dfn}

We prove the following proposition by using tensor nilpotence properties in the previous section.

\begin{prop}\label{main prop}
For $E,F\in{\rm DMF}(X,L,W)$, if ${\rm Supp}(E)\subseteq{\rm Supp}(F)$, then $E\in\langle F\rangle^{\otimes}$.

\begin{proof}
  Let $f:\mathcal{O}_X\rightarrow F^{\vee}\otimes F$
be the adjunction morphism in ${\rm DMF}(X,L,0)$ induced by the adjoint pair $(-)\otimes F\dashv \mathcal{H}om(F,-)$. Set $G:=C(f)[-1]$, and let $a:G\rightarrow\mathcal{O}_X$ be a morphism which completes the following distinguished triangle
$$G\xrightarrow{a}\mathcal{O}_X\xrightarrow{f}F^{\vee}\otimes F\rightarrow G[1].$$
Then, since $\langle F\rangle^{\otimes}$ is closed under  ${\rm DMF}(X,L,0)$-action, $E\otimes C(a)\cong(E\otimes F^{\vee})\otimes F\in \langle F\rangle^{\otimes}.$ We claim that for any $n>0$, $E\otimes C(a^{\otimes n})\in \langle F\rangle^{\otimes}$. Indeed, consider the following diagram

\[\xymatrix{
G^{\otimes n+1}\ar[rr]^{G\otimes (a^{\otimes n})}\ar@{=}[d]\ar@{}[drr]|\circlearrowright&& G\ar[rr]^{G\otimes b_n}\ar[d]_{a}&& G\otimes C(a^{\otimes n})\ar@{-->}[d]\ar[rr]&&G^{\otimes n+1}[1]\ar@{=}[d]\\
G^{\otimes n+1}\ar[rr]^{a^{\otimes n+1}}&&\mathcal{O}_X\ar[rr]^{b_{n+1}}\ar[d]&&C(a^{\otimes n+1})\ar@{-->}[d]\ar[rr]&&G^{\otimes n+1}\\
&&C(a)\ar@{=}[rr]\ar[d]&&C(a)\ar@{-->}[d]&&\\
&&G[1]\ar[rr]^{G\otimes b_n[1]}&&G\otimes C(a^{\otimes n})[1]&&
}\]
 where the top horizontal sequence is the distinguished triangle obtained by tensoring $G$ with the following
  triangle
 $$(\ast)\hspace{3mm}G^{\otimes n}\xrightarrow{a^{\otimes n}}\mathcal{O}_X\xrightarrow{b_n} C(a^{\otimes n})\rightarrow G^{\otimes n}[1].$$
 Then, by the octahedral axiom, we obtain the  triangle completing the vertical sequence on the right side in the above diagram 
 $$G\otimes C(a^{\otimes n})\rightarrow C(a^{\otimes n+1})\rightarrow C(a)\rightarrow G\otimes C(a^{\otimes n})[1].$$
Considering the triangle obtained by tensoring this triangle with $E$, we can prove that $E\otimes C(a^{\otimes n})\in \langle F\rangle^{\otimes}$ by induction on $n$. 

Tensoring $E$ with the above triangle $(\ast)$ for any $n>0$, we have a triangle
$$E\otimes G^{\otimes n}\xrightarrow{E\otimes a^{\otimes n}} E\longrightarrow  E\otimes C(a^{\otimes n})\longrightarrow E\otimes G^{\otimes n}[1].$$
If  $E\otimes a^{\otimes n}=0$ for some $n>0$, then $E$ is a direct summand of $E\otimes C(a^{\otimes n})\in \langle F\rangle^{\otimes}$, which implies that $E\in \langle F\rangle^{\otimes}$. Hence it suffices to show that there is an integer $n>0$ such that $E\otimes a^{\otimes n}=0$. By Lemma \ref{tens nilp2}, it is enough to show that for any $p\in {\rm Supp}(E)$,  $a\otimes k(p)=0$ in ${\rm KMF}(k(p),0)$. Since $p\in{\rm Supp}(E)\subseteq {\rm Supp}(F)$, $p\in {\rm Supp}(F^{\vee}\otimes F)$. By Lemma \ref{stalk lem2}, $F^{\vee}\otimes F\otimes k(p)\cong\mathcal{H}om_{k(p)}(F\otimes k(p),F\otimes k(p))\ne 0$ in ${\rm KMF}(k(p),0)$. Hence the natural map
$$k(p)\xrightarrow{g} \mathcal{H}om_{k(p)}(F\otimes k(p),F\otimes k(p))$$
is a split mono by Corollary \ref{stalk lem}, since it is non-zero map. Since $f\otimes k(p):k(p)\rightarrow  F^{\vee}\otimes F\otimes k(p)$ is equal to the composition of $g:k(p)\rightarrow\mathcal{H}om_{k(p)}(F\otimes k(p),F\otimes k(p))$ and the natural isomorphism $\mathcal{H}om_{k(p)}(F\otimes k(p),F\otimes k(p))\cong F^{\vee}\otimes F\otimes k(p)$, $f\otimes k(p)$ is also a split mono. Hence the triangle 
$$G\otimes k(p)\xrightarrow{a\otimes k(p)}k(p)\xrightarrow{f\otimes k(p)} F^{\vee}\otimes F\otimes k(p)\rightarrow G\otimes k(p)[1]$$
implies that $a\otimes k(p)=0$.
\end{proof}
\end{prop}

Now we are ready to prove the following main result. Recall that a subset $S\subseteq T$ of a topological space $T$ is called \textbf{specialization-closed} if it is a union of closed subsets of $T$. We easily see that $S$ is specialization-closed if and only if $s\in S$ implies $\overline{\{s\}}\subseteq S$. 

\begin{thm}\label{main result}
Let $X$ be a separated Noetherian scheme with an ample family of line bundles, $L$ be a line bundle on $X$, and $W\in \Gamma(X,L)$ be a non-zero-divisor.  There is one-to-one correspondence:
$$
\biggl\{
\begin{aligned}
\hspace{1mm}specialization{\normalbar}closed \hspace{3mm}\\subsets\hspace{1mm} of \hspace{1mm}\Sing(X_0/X)\hspace{1mm}
\end{aligned}
\biggr\}
\begin{aligned}
\xrightarrow{\sigma}\\
\xleftarrow{\tau}\\
\end{aligned}
\biggl\{
\begin{aligned}
\hspace{1mm}thick \hspace{1mm}\otimes{\normalbar}submodules\hspace{1mm}\\ of \hspace{1mm}\DMF(X,L,W)\hspace{3mm}
\end{aligned}
\biggr\}
$$
The bijective map $\sigma$ sends $Y$ to the thick subcategory consisting of matrix factorizations $F\in {\rm DMF}(X,L,W)$ with ${\rm Supp}(F)\subseteq Y$.
The inverse bijection $\tau$ sends  $\mathcal{T}$  to the union $\bigcup_{F\in\mathcal{T}}{\rm Supp}(F)$.
\begin{proof}
The map $\sigma$ is well defined since for any $E\in{\rm DMF}(X,L,0)$ and $F\in{\rm DMF}(X,L,W)$, we have ${\rm Supp}(E\otimes F)\subseteq{\rm Supp}(E)\cap{\rm Supp}(F)$ by Lemma \ref{supp tens}. The map $\tau$ is also well defined by Proposition \ref{supp local}.(2) and Proposition \ref{rel sing = supp}.(1).  

We show that $\sigma$ and $\tau$ are mutually inverse. Let $Y$ be a specialization-closed subset of ${\rm Sing}(X_0/X)$, and let $\mathcal{T}$ be a thick $\otimes$-submodule of ${\rm DMF}(X,L,W)$. By construction, we have $\tau(\sigma(Y))\subseteq Y$ and $\mathcal{T}\subseteq\sigma(\tau(\mathcal{T}))$. Hence it is enough to show the inclusions $Y\subseteq \tau(\sigma(Y))$ and $\sigma(\tau(\mathcal{T}))\subseteq\mathcal{T}$. 
 
 Since $\Sing(X_0/X)$ is a specialization-closed subset of $X$ by Proposition \ref{rel sing = supp}.(1), $Y$ is  specialization-closed in $X$. Hence $Y$ can be described as a union of closed subsets $Y_{\lambda}$ of $X$; $Y=\bigcup Y_{\lambda}$. 
 By Proposition \ref{support exist}.(2), for each $Y_{\lambda}$, there exists $F_{\lambda}\in {\rm DMF}(X,L,W)$ with ${\rm Supp}(F_{\lambda})=Y_{\lambda}$. Since $F_{\lambda}\in \sigma(Y)$, we have $Y_{\lambda}={\rm Supp}(F_{\lambda})\subseteq\tau(\sigma(Y))$. Hence $Y\subseteq \tau(\sigma(Y))$.

 To finish the proof, we show that $\sigma(\tau(\mathcal{T}))\subseteq\mathcal{T}$. Let $F\in\sigma(\tau(\mathcal{T}))$ be an object. Then, by construction, we have ${\rm Supp}(F)\subseteq\bigcup_{T_{\lambda}\in \mathcal{T}}{\rm Supp}(T_{\lambda})$. Then, as in the proof of \cite[Theorem 3.15]{thomason}, there is a finite set  $\{T_\lambda\}_{\lambda\in\Lambda}$ of objects in $\mathcal{T}$ such that ${\rm Supp}(F)\subseteq\bigcup_{\lambda\in\Lambda}{\rm Supp}(T_{\lambda})$. Hence ${\rm Supp}(F)\subseteq {\rm Supp}(\oplus_{\lambda\in\Lambda}T_\lambda)$. Since $\oplus_{\lambda\in\Lambda}T_\lambda\in \mathcal{T}$, it follows from Proposition \ref{main prop} that $F\in \mathcal{T}$.
 \end{proof}
\end{thm}


\begin{dfn} 
We say that an object $G\in\DMF(X,L,W)$ is a \textbf{$\otimes$-generator} of $\DMF(X,L,W)$ if  $\langle G \rangle^{\otimes}=\DMF(X,L,W)$. 
\end{dfn}

The following corollary says that the closedness of  relative singular locus $\Sing(X_0/X)$ in $X_0$ is related to the existence of a $\otimes$-generator of $\DMF(X,L,W)$.

\begin{cor}\label{generator}
Notation is same as in Theorem \ref{main result}. The relative singular locus $\Sing(X_0/X)$ is closed in $X_0$ if and only if  $\DMF(X,L,W)$ has a $\otimes$-generator.
\begin{proof}
Assume that the subset  $\Sing(X_0/X)$ is closed in $X_0$. Since the relative singular locus $\Sing(X_0/X)$ is the union of supports of all objects in $\DMF(X,L,W)$ and $X$ is Noetherian, there is  a finite subset $\{F_i\}$ of objects $F_i\in\DMF(X,L,W)$ such that $\Sing(X_0/X)=\bigcup_{i}\Supp(F_i)=\Supp(\oplus_i F_i)$. By Theorem \ref{main result}, there is a specialization-closed subset $Y\subseteq \Sing(X_0/X)$ such that $\sigma(Y)=\langle\oplus_{i}F_i\rangle^{\otimes}$. Since $\Sing(X_0/X)=\Supp(\oplus_iF_i)\subseteq Y$, we have $Y=\Sing(X_0/X)$. Hence $\langle\oplus_{i}F_i\rangle^{\otimes}=\sigma(\Sing(X_0/X))=\DMF(X,L,W)$.

If $\DMF(X,L,W)$ has a $\otimes$-generator $G$,   for any object $F\in\DMF(X,L,W)$, we have $\Supp(F)\subseteq\Supp(G)$ by Lemma \ref{supp data}. Hence Proposition \ref{rel sing = supp}.(1) implies that  $\Sing(X_0/X)=\Supp(G)$, and so $\Sing(X_0/X)$ is closed in $X_0$.
 \end{proof}
\end{cor}




\begin{rem} Let $(X,L,W)$ be the same LG model as in Theorem \ref{main result}. If $X$ is regular and $L$ is ample, any thick subcategory of $\DMF(X,L,W)$ is automatically $\otimes$-submodule. In particular, the set  on the right-hand side in Theorem \ref{main result} is equal to the set of thick subcategories of $\DMF(X,L,W)$. Since this fact is proved by Stevenson  in a different context \cite{stev},  we do not  include the proof here.

\end{rem}

\begin{exa} 
Let $(X,W)$ and $X_0$ be same as in the Example \ref{example}.(2). Then  $\Sing(X_0/X)=\{ \langle \overline{x}, \overline{y}, \overline{z}-a\rangle \mid a\in\bC\setminus \{0\}\}$. In this case, using the above results we see the following:\vspace{1mm}\\
(1) Since any point in $\Sing(X_0/X)$ is a closed point,  by Theorem \ref{main result} the set of thick $\otimes$-submodules of $\DMF(X,W)$ is bijective to the set of arbitrary subsets of $\bC\setminus \{0\}$.\vspace{1mm}\\
(2) Since $\Sing(X_0/X)$ is  not a closed subset of $X_0$, by Corollary \ref{generator} $\DMF(X,W)$ does not have a $\otimes$-generator. In particular, $\DMF(X,W)$  does not have a classical generator, i.e. there is no object $F\in\DMF(X,W)$ such that $\langle F\rangle=\DMF(X,W)$. This was proved in \cite[Section 3.3]{efi-posi} by a different argument.

\end{exa}

\vspace{1mm}
\section{Tensor structures on matrix factorizations and its spectrum}

Using the classification result in the previous section, we will construct the relative singular loci from derived matrix factorization categories by considering tensor structures induced by tensor products.

 \subsection{Balmer's tensor triangular geometry}
Following \cite{balmer} and \cite[Chapter 4]{yu}, we will recall some basic definitions and results of the theory of tensor triangular geometry.

\begin{dfn} A \textbf{pseudo tensor triangulated category}  $(\mathcal{T},\otimes)$ consists of  a triangulated category $\mathcal{T}$ and  symmetric associated bifunctor $\otimes: \mathcal{T}\times\mathcal{T}\rightarrow\mathcal{T}$ which is exact in each variable. For the precise definition, see \cite[Definition 4.1.1]{yu}. 
\end{dfn}

\begin{rem}
We don't assume that a pseudo tensor triangulated category has  a unit $1_{\otimes}$,  and this  is the only difference from the original definition of tensor triangulated categories in \cite{balmer}.
\end{rem}

\begin{dfn}
Let $(\mathcal{T},\otimes)$ be a pseudo tensor triangulated category. 
\begin{itemize}
\item[$(1)$] A thick subcategory $\mathcal{I}\subset \mathcal{T}$ is called \textbf{$\otimes$-ideal} if the following implication holds:
\begin{center}
$A\in\mathcal{T}$ and $B\in\mathcal{I}$ \,\,$\Rightarrow$ $A\otimes B\in\mathcal{I}$.
\end{center}
\item[$(2)$] A $\otimes$-ideal $\mathcal{P}$ is called \textbf{prime} if the following holds
\begin{center}
$A\notin\mathcal{P}$ and $B\notin\mathcal{P}$ \,\,$\Rightarrow$ $A\otimes B\notin\mathcal{P}$.
\end{center}
\item[$(3)$] A $\otimes$-ideal $\mathcal{I}$ is called \textbf{radical} if $\sqrt{\mathcal{I}}=\mathcal{I}$, where $\sqrt{\mathcal{I}}$ is the radical of $\mathcal{I}$, i.e. 
$$\sqrt{\mathcal{I}}:=\{\,A\in\mathcal{T}\,\,|\,\exists\, n\geq1\,\,\,{\rm such \,\,that\,\,} A^{\otimes n}\in\mathcal{I}\,\}$$
\end{itemize}
\end{dfn}

\vspace{1mm}
For a pseudo tensor triangulated category $(\mathcal{T},\otimes)$, we can consider  the Zariski topology on the set of all prime ideals of $(\mathcal{T},\otimes)$.

\begin{dfn}The \textbf{spectrum}, denoted by ${\rm Spc}(\mathcal{T},\otimes)$, of $(\mathcal{T},\otimes)$ is defined as the set of all prime $\otimes$-ideals
\begin{center}
${\rm Spc}(\mathcal{T},\otimes):=\{\,\mathcal{P}\,|\, \mathcal{P}$ is a prime $\otimes$-ideal $\}$\vspace{1mm}
\end{center}
The \textbf{Zariski topology} on ${\rm Spc}(\mathcal{T},\otimes)$ is defined by the collection of closed subsets of the form $Z(\mathcal{S}):=\{\mathcal{P}\in{\rm Spc}(\mathcal{T},\otimes)\,|\,\mathcal{S}\cap\mathcal{P}=\emptyset\}$ for any family of objects $\mathcal{S}\subseteq\mathcal{T}$.
\end{dfn}

\begin{dfn}\label{df of supp data}
A \textbf{support data} on a pseudo tensor triangulated category $(\mathcal{T}, \otimes)$ is  a  pair $(X,\sigma)$ of a topological space $X$ and an assignment  $\sigma:{\rm Ob}\mathcal{T}\rightarrow\{{\rm closed\hspace{1mm}subsets\hspace{1mm} of \hspace{1mm}} X \}$ satisfying the following conditions:
\begin{itemize}
\item[$(1)$] $\sigma(0)=\emptyset$ and $\bigcup_{A\in\mathcal{T}}\sigma(A)=X$.
\item[$(2)$] $\sigma(A\oplus B)=\sigma(A)\cup\sigma(B)$.
\item[$(3)$] $\sigma(A[1])=\sigma(A)$.
\item[$(4)$] $\sigma(A)\subseteq\sigma(B)\cup\sigma(C)$ for any triangle $A\rightarrow B\rightarrow C\rightarrow A[1]$.
\item[$(5)$] $\sigma(A\otimes B)=\sigma(A)\cap\sigma(B).$
\end{itemize}
We say that a support data $(X,\sigma)$ is \textbf{classifying} if the following properties hold:
\begin{itemize}
\item[$({\rm a})$] The topological space $X$ is Noetherian and any non-empty irreducible closed subset has a unique generic point. 
\item[$({\rm b})$] There is the following bijective correspondence:
$$\Theta:\bigl\{ \, {\rm specialization{\normalbar}closed \,\,subsets\,\,of\,\,}X\bigr\}\xrightarrow{\sim}\bigl\{\,{\rm radical \,\,thick \,\,}\otimes\mathchar`-{\rm ideals}\,\bigr\}$$
defined by $\Theta(Y):=\{A\in\mathcal{T}\,|\,\sigma(A)\subseteq X\}$, with $\Theta^{-1}(\mathcal{I})=\bigcup_{A\in\mathcal{I}}\sigma(A)$.
\end{itemize}
\end{dfn}

\begin{rem}
Because of the lack of the unit $1_{\otimes}$ in $(\mathcal{T}, \otimes)$ in our setting, we replace the condition $\sigma(1_{\otimes})=X$ in the original definition of  support data in \cite[Definition 3.1 (SD1)]{balmer} with $\bigcup_{A\in\mathcal{T}}\sigma(A)=X$.

\end{rem}

The following result is essentially due to Balmer \cite{balmer}, and it is the key result for the result in the next subsection. See also \cite[Theorem 4.1.16]{yu}.

\begin{thm}[{\cite[Theorem 5.2]{balmer}}]\label{homeo}
Assume that $(X,\sigma)$ is a classifying support data on a pseudo tensor triangulated category $(\mathcal{T},\otimes)$. Then we have the canonical homeomorphism
$$f:X\xrightarrow{\sim}{\rm Spc}((\mathcal{T},\otimes)),$$
defined by $f(x):=\{A\in\mathcal{T}\,|\,x\notin\sigma(A)\}$.
\end{thm}

\subsection{Construction of relative singular loci from matrix factorizations}

In this section, using our classification result, we construct relative singular loci from pseudo tensor triangulated structures on derived matrix factorization categories. This kind of observation is also discussed in \cite{yup}.
Following \cite{yup}, we consider the natural pseudo tensor triangulated structure on derived matrix factorization categories. 

Throughout this section,  $X$ is a separated Noetherian scheme with an ample family of line bundles, $L$ is a line bundle on $X$, and $W\in\Gamma(X,L)$ is a non zero-divisor. Denote by $X_0$ the zero scheme of $W$, and assume that $2\in \Gamma(X,\mathcal{O}_X)$ is a unit of the ring $\Gamma(X,\mathcal{O}_X)$.

\vspace{2mm}
For any unit $\lambda\in\Gamma(X,\mathcal{O}_X)^{\times}$ in the ring $\Gamma(X,\mathcal{O}_X)$, we have  a natural functor $$\lambda:{\rm DMF}(X,L,W)\rightarrow{\rm DMF}(X,L,\lambda W)$$ 
defined by $\lambda\bigl(F_1\xrightarrow{\varphi_1}F_0\xrightarrow{\varphi_0}F_1\bigr):=\bigl(F_1\xrightarrow{\varphi_1}F_0\xrightarrow{\lambda\varphi_0}F_1\bigr)$, and it is a triangulated equivalence (see \cite[Proposition 4.1.19]{yu}).

\begin{dfn}
Suppose that $2\in \Gamma(X,\mathcal{O}_X)$ is a unit. Then we define a bifunctor 
$$\otimes^{\frac{1}{2}}:{\rm DMF}(X,L,W)\times{\rm DMF}(X,L,W)\rightarrow{\rm DMF}(X,L,W)$$
as the composition ${\rm DMF}(X,L,W)\times{\rm DMF}(X,L,W)\xrightarrow{\otimes}{\rm DMF}(X,L,2W)\xrightarrow{\frac{1}{2}}{\rm DMF}(X,L,W)$.

\end{dfn}

As in \cite[Proposition 4.1.22]{yu}, we see that the pair $({\rm DMF}(X,L,W),\otimes^{\frac{1}{2}})$ is a pseudo tensor triangulated category.

Consider an assignment 
\[{\rm Supp(-)}:{\rm Ob}({\rm DMF}(X,L,W))\rightarrow\{{\rm specialization{\normalbar}closed \,\,subsets\,\, of\,\,} {\rm Sing}(X_0/X)\}\]
 defined by the supports of objects in $\DMF(X,L,W)$. 

\begin{thm}\label{cl supp data}
$\bigl({\rm Supp}, \,{\rm Sing}(X_0/X)\bigr)$ is a classifying support data on $\bigl({\rm DMF}(X,L,W),\otimes^{\frac{1}{2}}\bigr)$.
\begin{proof}
Since the functor $\frac{1}{2}:{\rm DMF}(X,L,2W)\rightarrow{\rm DMF}(X,L,W)$ is an equivalence and commutes with taking stalks, we have 
$${\rm Supp}(\frac{1}{2}(F))={\rm Supp}(F).$$ 
Therefore, it follows from Lemma \ref{supp data} and Proposition \ref{rel sing = supp} that $\bigl({\rm Supp}, {\rm Sing}(X_0/X)\bigr)$ is a support data. We will show that the support data satisfies the conditions (a) and (b) in Definition \ref{df of supp data}.

We check the condition (a). Note that $X$ satisfies the condition (a). It follows that  ${\rm Sing}(X_0/X)$ is a Noetherian topological space, since so is $X$. Note that any {\it irreducible} closed subset $Z$ of $\Sing(X_0/X)$ is closed in $X$. Indeed, since $Z$ is closed in a specialization-closed subset of $X$, $Z$ is specialization-closed in $X$. Hence $Z$ is a union of irreducible closed subsets $Z_{\lambda}$ of $X$; $Z=\bigcup_{\lambda} Z_{\lambda}$. Since $Z$ is irreducible in $\Sing(X_0/X)$, there is an irreducible closed subset $Z_{\lambda'}$ of $X$ such that $Z=Z_{\lambda'}$. Hence $Z$ is an irreducible closed subset of $X$, and it has a unique generic point.

Next, we verify the condition (b). By Theorem \ref{main result},
 it is enough to show that for a thick subcategory $\mathcal{T}$ of ${\rm DMF}(X,L,W)$
\begin{center}
$\mathcal{T}$ is $\otimes$-submodule $\Leftrightarrow$ $\mathcal{T}$ is radical $\otimes^{\frac{1}{2}}$-ideal.
\end{center}
The implication $(\Rightarrow)$ follows immediately from Theorem \ref{main result} and Lemma \ref{supp data}.(4).
We show the other implication $(\Leftarrow)$. For this, let $\mathcal{I}\subset \DMF(X,L,W)$ be a radical thick $\otimes^{\frac{1}{2}}$-ideal.  For objects $E\in{\rm DMF}(X,L,0)$ and $F\in\mathcal{I}$, it suffices  to prove that $E\otimes F\in\mathcal{I}$. We have 
\begin{align}
(E\otimes F)\otimes^{\frac{1}{2}}(E\otimes F)&=\frac{1}{2}\bigl((E\otimes F)\otimes(E\otimes F)\bigr)\notag\\
&\cong \frac{1}{2}\bigl((E\otimes F\otimes E)\otimes F)\notag\\
&=(E\otimes F\otimes E)\otimes^{\frac{1}{2}} F,\notag
\end{align}
and the object in the bottom line is in $\mathcal{I}$ since $E\otimes F\otimes E\in{\rm DMF}(X,L,W)$ and $\mathcal{I}$ is $\otimes^{\frac{1}{2}}$-ideal. Hence  $E\otimes F\in\mathcal{I}$ as $\mathcal{I}$ is radical.
\end{proof}
\end{thm}

Theorem \ref{homeo} and Theorem \ref{cl supp data} imply the following result.

\begin{cor}\label{reconst}
There is a homeomorphism
$${\rm Spc}({\rm DMF}(X,L,W),\otimes^{\frac{1}{2}})\cong{\rm Sing}(X_0/X)$$

\end{cor}

\begin{rem}
By Proposition \ref{affine} and Proposition \ref{rel sing = sing}, we see that Corollary \ref{reconst} is a generalization of \cite[Theorem 1.2]{yup}, where Yu consider the case when $X$ is a  regular affine scheme of finite Krull dimension. 
\end{rem}

\vspace{4mm}
\address{Department of Mathematics, Kyoto University, Kitashirakawa-Oiwake-cho, Sakyo-ku, Kyoto, 606-8502, Japan}\\
{\it E-mail address}: \email{y.hirano@math.kyoto-u.ac.jp}

\end{document}